\documentclass[12pt,twoside,a4paper,openany]{article}
\usepackage[centertags]{amsmath}
\usepackage{amsfonts}
\usepackage[all]{xy}
\usepackage{graphicx}
\usepackage{amsthm}
\usepackage{newlfont}
\usepackage[latin1]{inputenc}
\usepackage{amsmath}
\usepackage{eucal}
\usepackage{amssymb}
\usepackage{titlesec}
\titleformat{\subsection}[runin]
       {\bfseries}
			 {\thesubsection}
       {0.5em}
       {}
       []

\newtheorem{thm}{Theorem}[section]
\newtheorem{cor}[thm]{Corollary}

\newtheorem{prop}[thm]{Proposition}

\newtheorem{subs}[thm]{}
\newtheorem{rem}[thm]{Remark}

\newtheorem{defn}[thm]{Definition}

\begin{document}

\begin{center}
\textbf{\Large Simple limit linear series for curves of compact type}

\textbf{\Large with three irreducible components}
\end{center}

\begin{center}
\text{Author: Gabriel Armando Mu\~noz Márquez}
\end{center}

\begin{center}
\text{Affiliations:}\\ 
\text{Universidad Científica del Sur, Lima, Perú.}\\
\text{Universidad Nacional Mayor de San Marcos, Lima, Perú.}
\end{center}

\begin{center}
\text{email: gmunozm@cientifica.edu.pe}
\end{center}

\begin{abstract}

We study simple Osserman limit linear series (that is, Osserman limit linear series having a simple basis) on curves of compact type with three irreducible components. For compact type curves with two components, every exact limit linear series is simple. But, for the case of three components, this property is no longer true. We study a certain distributivity property related with the existence of a simple basis and we find characterizations of that property at any fixed multidegree. Also, we find a characterization of simple limit linear series among the exact limit linear series. In our approach, we also study other structure similar to a simple basis and we get a certain inequality for any exact limit linear series. 

\end{abstract}

\textbf{Keywords} Compact type curves; simple limit linear series 

\textbf{Mathematics Subject Classification} 14H50

\section{Introduction}

For curves of compact type, Eisenbud and Harris [3] developed the theory of limit linear series as an analogue of linear series. This theory is very powerful for degeneration arguments on curves. The idea is to analyze how linear series degenerate when a family of smooth curves degenerates to a compact type curve. Eisenbud and Harris approached this situation by considering only the possible limit line bundles with nonnegative multidegree and degree $d$ on one irreducible component of the curve. Osserman [7] developed a new and more functorial construction for the theory of limit linear series. The basic idea is to consider all possible limit line bundles with nonnegative multidegree. Thus, Osserman limit linear series carry more information about limit line bundles. This new theory has a generalization to higher rank vector bundles Osserman [9].

Esteves and Osserman [4] related limit linear series to fibers of Abel maps via the definition of limit linear series by Osserman [7]. They restricted their work to the case of curves of compact type with two irreducible components and studied exact limit linear series (these contain in particular all limits of linear series on the generic fiber in a regular smoothing family). For any exact limit linear series, they establihed the existence of a simple basis (a certain structure of multidegrees and linearly independent sections), which was a key property in their work. For instance, the simple basis was useful to compute the Hilbert polynomial of certain closed subscheme of a fiber of an Abel map.

For curves of compact type with three irreducible components, Muñoz [5] generalized some of the results in Esteves and Osserman [4] for the special case of exact limit linear series arising as the unique exact extension of a refined Eisenbud-Harris limit linear series. Also, for any exact limit linear series satisfying a certain distributivity property, Muñoz [5] constructed a certain structure similar to a simple basis (Muñoz [5], Proposition 4.4). This structure played a similar role to that of the simple basis in Esteves and Osserman [4], so such a structure was enough for that work. However, we expected that such a structure could be used to construct a simple basis. In the present work, we do that construction (see Proposition \ref{simple basis existence}). In a different context, Dos Santos [2] defines a similar property to the distributivity property (there called the intersection property). He also related that property with the existence of a simple basis.

Our aim in this work is to study simple limit linear series (that is, limit linear series having a simple basis) on curves of compact type with three irreducible components. We find characterizations of the distributivity property at any fixed multidegree. Also, we find a characterization of simple limit linear series among the exact limit linear series. In our approach, we also get a certain inequality for any exact limit linear series.

We now explain the contents of the paper. We begin with the notation. In Esteves and Osserman [4], a limit linear series of degree $d$ and dimension $r$ on a curve of compact type $X$ with two irreducible components $Y$ and $Z$, meeting transversally at a point $P$, is a collection $\mathfrak{g}:=(\mathcal{L},V_{0},\ldots,V_{d})$, where $\mathcal{L}$ is an invertible sheaf on $X$ of degree $d$ on $Y$ and degree $0$ on $Z$, and $V_{i}$ is a vector subspace of $H^{0}(X,\mathcal{L}^{i})$ of dimension $r+1$, for each $i=0,\ldots,d$, where $\mathcal{L}^{i}$ is the invertible sheaf on $X$ with restrictions $\mathcal{L}\big|_{Y}(-iP)$ and $\mathcal{L}\big|_{Z}(iP)$, and these vector subspaces are linked by certain natural maps between the sheaves $\mathcal{L}^{i}$. Thus, a limit linear series is defined by a collection of pairs $(\mathcal{L}^{i}, V_{i})$, for each $i=0,\ldots,d$, and we can use the notation $\{(\mathcal{L}^{i}, V_{i})\}_{i}$, where $0\leq i\leq d$. For each $i=0,\ldots,d$, the invertible sheaf $\mathcal{L}^{i}$ has multidegree $\underline{d}:=(d-i,i)$. So, for each $i=0,\ldots,d$, setting $\underline{d}:=(d-i,i)$, $\mathcal{L}_{\underline{d}}=\mathcal{L}^{i}$ and $V_{\underline{d}}=V_{i}$, a limit linear series can also be denoted by a collection $\{(\mathcal{L}_{\underline{d}},V_{\underline{d}})\}_{\underline{d}}$, where $\underline{d}\geq \underline{0}$ of total degree $d$. We will use this notation for the case of three irreducible components. 

In this work, $X$ will denote the union of three smooth curves $X_1, X_2$ and $X_3$, such that $X_1$ and $X_2$ meet transversally at a point $A$, and $X_2$ and $X_3$ meet transversally at a point $B$, with $A\neq B$. A limit linear series of degree $d$ and dimension $r$ on $X$  is a collection $\mathfrak{g}:=\{(\mathcal{L}_{\underline{d}},V_{\underline{d}})\}_{\underline{d}}$, where $\underline{d}\geq \underline{0}$ of total degree $d$, each $\mathcal{L}_{\underline{d}}$ is an invertible sheaf of multidegree $\underline{d}$, and $V_{\underline{d}}$ is a vector subspace of $H^{0}(X,\mathcal{L}_{\underline{d}})$ of dimension $r+1$, for each $\underline{d}$, where $\mathcal{L}_{\underline{d}}$ is the invertible sheaf on $X$ with restrictions $\mathcal{L}\big|_{X_1}(-(d-i)A)$, $\mathcal{L}\big|_{X_2}((d-i)A-lB)$ and $\mathcal{L}\big|_{X_3}(lB)$, with $\underline{d}=(i,d-i-l,l)$ and $\mathcal{L}:=\mathcal{L}_{(d,0,0)}$, and the subspaces $V_{\underline{d}}$ are linked by certain maps between the sheaves $\mathcal{L}_{\underline{d}}$ (see Subsection \ref{limit1}). 

If $Y$ is a subcurve of $X$, for any $\underline{d}$ and for any subspace $V\subseteq H^{0}(X,\mathcal{L}_{\underline{d}})$, we denote by $V^{Y,0}$ the subspace of $V$ of sections that vanish on $Y$.

In Section 3, we define the distributivity property. In Proposition \ref{distributivity characterization}, which is one of our main results, we establish chracterizations of distributivity at a fixed multidegree under specific exactness conditions. Specifically, our Proposition \ref{distributivity characterization} says:

{\em For any $i\geq 0$ and $l\geq 0$ such that $i+l\leq d$, let $\underline{d}:=(i,d-i-l,l)$. Then, the following statements hold:\\
1. If $i>0$, let $\underline{d}'':=(i-1,d-i-l+1,l)$. Let $V_{\underline{d}}$ and $V_{\underline{d}''}$ be $r+1$-dimensional subspaces of $H^{0}(\mathcal{L}_{\underline{d}})$ and $H^{0}(\mathcal{L}_{\underline{d}''})$, respectively, such that $\varphi_{\underline{d},\underline{d}''}(V_{\underline{d}})=V_{\underline{d}''}^{X_{1},0}$. Then the distributivity holds in $V_{\underline{d}''}$ if and only if}
\begin{center}
$\text{dim}\,(V_{\underline{d}}^{X_{2},0}+V_{\underline{d}}^{X_{3},0})-\text{dim}\,(V_{\underline{d}''}^{X_{2},0}+V_{\underline{d}''}^{X_{3},0})=r+1-\text{dim}\,(V_{\underline{d}''}^{X_{1},0}+V_{\underline{d}''}^{X_{2},0}+V_{\underline{d}''}^{X_{3},0})$
\end{center}
{\em 2. If $i+l\leq d-1$, let $\underline{d}''':=(i,d-i-l-1,l+1)$. Let $V_{\underline{d}}$ and $V_{\underline{d}'''}$ be $r+1$-dimensional subspaces of $H^{0}(\mathcal{L}_{\underline{d}})$ and $H^{0}(\mathcal{L}_{\underline{d}'''})$, respectively, such that $\varphi_{\underline{d}''',\underline{d}}(V_{\underline{d}'''})=V_{\underline{d}}^{X_{3},0}$. Then the distributivity holds in $V_{\underline{d}}$ if and only if}
\begin{center}
$\text{dim}\,(V_{\underline{d}'''}^{X_{1},0}+V_{\underline{d}'''}^{X_{2},0})-\text{dim}\,(V_{\underline{d}}^{X_{1},0}+V_{\underline{d}}^{X_{2},0})=r+1-\text{dim}\,(V_{\underline{d}}^{X_{1},0}+V_{\underline{d}}^{X_{2},0}+V_{\underline{d}}^{X_{3},0})$,
\end{center}
where the maps $\varphi_{\underline{d},\underline{d}''}$ and $\varphi_{\underline{d}''',\underline{d}}$ in Proposition \ref{distributivity characterization} are the maps linking the corresponding sheaves. Another important result is Proposition \ref{series exactas desigualdad}, which establishes an inequality for any exact limit linear series. Specifically, our Proposition \ref{series exactas desigualdad} says:

{\em Let $\{(\mathcal{L}_{\underline{d}},V_{\underline{d}})\}_{\underline{d}}$ be an exact limit linear series of degree $d$ and dimension $r$. Then} 
\begin{center}
$\displaystyle\sum_{\underline{d}}\text{dim}\,\left(\dfrac{V_{\underline{d}}}{V_{\underline{d}}^{X_{1},0}+V_{\underline{d}}^{X_{2},0}+V_{\underline{d}}^{X_{3},0}}\right)\geq r+1$.
\end{center}
As a consequence of Proposition \ref{distributivity characterization} and Proposition \ref{series exactas desigualdad}, in Corollary \ref{sum of codimensions}, we get the following characterization of exact limit linear series satisfying the distributivity at each multidegree:\\
{\em Let $\{(\mathcal{L}_{\underline{d}},V_{\underline{d}})\}_{\underline{d}}$ be an exact limit linear series of degree $d$ and dimension $r$. Then} 
\begin{center}
$\displaystyle\sum_{\underline{d}}\text{dim}\,\left(\dfrac{V_{\underline{d}}}{V_{\underline{d}}^{X_{1},0}+V_{\underline{d}}^{X_{2},0}+V_{\underline{d}}^{X_{3},0}}\right)=r+1$\\
{\em if and only if the distributivity holds in} $V_{\underline{d}}$ {\em for any} $\underline{d}$.
\end{center}
In Section 3, we also study simple limit linear series. In Proposition \ref{multidegrees characterization}, we characterize the multidegrees which determine all simple bases. In Proposition \ref{like simple basis}, following the construction in Muñoz [5], Proposition 4.4, we present a certain structure of linearly independent sections for any exact limit linear series satisfying the distributivity property. In propositions \ref{prop1structure} and \ref{prop2structure}, we obtain some properties of that structure. These properties imply the existence of a simple basis (we give a proof of that in Proposition \ref{simple basis existence}).

We conclude Section 3 with a characterization of simple limit linear series among the exact limit linear series. Specifically, our Corollary \ref{simple characterization} says:\\ 
{\em Let} $\{(\mathcal{L}_{\underline{d}},V_{\underline{d}})\}_{\underline{d}}$ {\em be an exact limit linear series of degree} $d$ {\em and dimension} $r$. {\em Then the following statements are equivalent:}\\
1. $\{(\mathcal{L}_{\underline{d}},V_{\underline{d}})\}_{\underline{d}}$ {\em is simple}.\\
2. $\displaystyle\sum_{\underline{d}}\text{dim}\,\left(\dfrac{V_{\underline{d}}}{V_{\underline{d}}^{X_{1},0}+V_{\underline{d}}^{X_{2},0}+V_{\underline{d}}^{X_{3},0}}\right)=r+1$.\\
3. {\em The distributivity holds in} $V_{\underline{d}}$ {\em for any} $\underline{d}$.

In a subsequent work, given a refined limit linear series, we will construct all simple extensions by combining our Proposition \ref{distributivity characterization} and the method of Muñoz [6] for the construction of all exact extensions.

\section{Preliminaries}\label{limit}
\begin{subs}{\em(}Limit linear series{\em)}\label{limit1} {\em Throughout this article, $X$ will denote the union of three smooth curves $X_1, X_2$ and $X_3$, such that $X_1$ and $X_2$ meet transversally at a point $A$, and $X_2$ and $X_3$ meet transversally at a point $B$, with $A\neq B$. If $Y$ is a reduced union of some components of $X$, we denote by $Y^{c}$ the closure of the complement of Y. Also, note that if $Y$ is a reduced union of some components of $X$, we get the following exact sequence
\[0\rightarrow \mathcal{L}\big|_{Y^{c}}(-Y\cap Y^{c})\rightarrow \mathcal{L}\rightarrow \mathcal{L}\big|_{Y}\rightarrow 0,\]
for any invertible sheaf $\mathcal{L}$ on $X$.

Notice that, since $X$ is a curve of compact type, invertible sheaves on $X$ are uniquely determined up to isomorphism by their restrictions to $X_1, X_2$ and $X_3$. Fix an invertible sheaf $\mathcal{O}_{(X_1)}$ on $X$ with restrictions $\mathcal{O}_{X_1}(-A)$, $\mathcal{O}_{X_2}(A)$ and $\mathcal{O}_{X_3}$. Analogously, fix an invertible sheaf $\mathcal{O}_{(X_3)}$ on $X$ with restrictions $\mathcal{O}_{X_1}$, $\mathcal{O}_{X_2}(B)$ and $\mathcal{O}_{X_3}(-B)$. We also fix an invertible sheaf $\mathcal{O}_{(X_{2})}$ on $X$ with restrictions $\mathcal{O}_{X_1}(A)$, $\mathcal{O}_{X_2}(-A-B)$ and $\mathcal{O}_{X_3}(B)$. We observe that $\mathcal{O}_{(X_1)}\otimes \mathcal{O}_{(X_2)}\otimes \mathcal{O}_{(X_3)}\cong \mathcal{O}_{X}$. Fix such an isomorphism, and call $\tau$ this map. For each $q=1,2,3$, set
\[\mathcal{O}_{(-X_{q})}:=\bigotimes_{\tilde{q}\neq q} \mathcal{O}_{(X_{\tilde{q}})}.\]
Given an invertible sheaf $\mathcal{L}$ on $X$ of degree $d$ on $X_1$ and degree $0$ on $X_2$ and $X_3$,
for each $i\geq 0$ and $l\geq 0$ such that $i+l\leq d$, set 
\begin{center}
$\mathcal{L}_{(i,d-i-l,l)}:=\mathcal{L}\otimes \mathcal{O}_{(X_1)}^{\otimes d-i}\otimes \mathcal{O}_{(-X_3)}^{\otimes l}$. 
\end{center}
We observe that $\mathcal{L}_{(i,d-i-l,l)}$ has multidegree $(i,d-i-l,l)$. 
For any $i\geq 0$ and $l\geq 0$ such that $i+l\leq d$, let $\underline{d}:=(i,d-i-l,l)$. For each $q=1,2,3$, set
 \[\widetilde{\underline{d}}_{q}:=\left\{
            \begin{array}{rcl} (i-1,d-i-l+1,l) & \mbox{if} & q=1,\\
               (i+1,d-i-l-2,l+1) & \mbox{if} & q=2,\\
							 (i,d-i-l+1,l-1) & \mbox{if} & q=3.
                           \end{array}\right.\]													
Thus, $\widetilde{\underline{d}}_{q}=\underline{d}+\text{multideg}(\mathcal{O}_{(X_{q})})$ for each $q=1,2,3$.

For each $\underline{d}\geq \underline{0}$ and $\widetilde{\underline{d}}_{q}\geq \underline{0}$, there exist maps $\varphi_{\underline{d},\widetilde{\underline{d}}_q}$ and $\varphi_{\widetilde{\underline{d}}_q,\underline{d}}$ (see Muñoz [6]) satisfying the following properties:\\
(1) $\varphi_{\underline{d},\widetilde{\underline{d}}_q}$ factors as
\[\varphi_{\underline{d},\widetilde{\underline{d}}_q}: \mathcal{L}_{\underline{d}}\rightarrow \mathcal{L}_{\underline{d}}\big|_{X_{q}^{c}}\cong\mathcal{L}_{\widetilde{\underline{d}}_q}\big|_{X_{q}^{c}}(-X_{q}\cap X_{q}^{c})\hookrightarrow \mathcal{L}_{\widetilde{\underline{d}}_{q}},\]
where the first map in the composition is the restriction map and the last map is the natural inclusion,\\
(2) $\varphi_{\widetilde{\underline{d}}_q,\underline{d}}$ factors as
\[\varphi_{\widetilde{\underline{d}}_q,\underline{d}}: \mathcal{L}_{\widetilde{\underline{d}}_q}\rightarrow \mathcal{L}_{\widetilde{\underline{d}}_q}\big|_{X_{q}}\cong\mathcal{L}_{\underline{d}}\big|_{X_{q}}(-X_{q}\cap X_{q}^{c})\hookrightarrow \mathcal{L}_{\underline{d}},\]
where the first map in the composition is the restriction map and the last map is the natural inclusion, and\\
(3)
\begin{center}
$\varphi_{\underline{d}^{2},\underline{d}^{3}}\circ \varphi_{\underline{d}^{1},\underline{d}^{2}}=\varphi_{\underline{d}^{1},\underline{d}^{3}}$,
\end{center}
for $\underline{d}^{2}=\underline{d}^{1}+\text{multideg}(\mathcal{O}_{(X_{q_1})})$ and $\underline{d}^{3}=\underline{d}^{2}+\text{multideg}(\mathcal{O}_{(X_{q_2})})$ with $q_{1}\neq q_{2}$.

We observe that the compositions $\varphi_{\underline{d},\widetilde{\underline{d}}_q}\circ \varphi_{\widetilde{\underline{d}}_q,\underline{d}}$ and $\varphi_{\widetilde{\underline{d}}_q,\underline{d}}\circ \varphi_{\underline{d},\widetilde{\underline{d}}_q}$ are zero.

If $Y$ is a subcurve of $X$, for any $\underline{d}$ and for any subspace $V\subseteq H^{0}(X,\mathcal{L}_{\underline{d}})$, we denote by $V^{Y,0}$ the subspace of $V$ of sections that vanish on $Y$. If $Y$ is an irreducible component of $X$, we denote by $\mathcal{L}_{Y}$ the invertible sheaf $\mathcal{L}_{\underline{d}}$, where the component of $\underline{d}$ corresponding to $Y$ is equal to $d$ and the other components of $\underline{d}$ are $0$. Also, to ease notation, let $\mathcal{L}_{il}:=\mathcal{L}_{(i,d-i-l,l)}$.

Fix integers $d$ and $r$. A {\em limit linear series} on $X$ of degree $d$ an dimension $r$ is a collection consisting of an invertible sheaf $\mathcal{L}$ on $X$ of degree $d$ on $X_1$ and degree $0$ on $X_2$ and $X_3$, and vector subspaces $V_{\underline{d}}\subseteq H^{0}(X,\mathcal{L}_{\underline{d}})$ of dimension $r+1$, for each $\underline{d}:=(i,d-i-l,l)\geq \underline{0}$, such that $\varphi_{\underline{d},\widetilde{\underline{d}}_q}(V_{\underline{d}})\subseteq V_{\widetilde{\underline{d}}_q}$ and $\varphi_{\widetilde{\underline{d}}_q,\underline{d}}(V_{\widetilde{\underline{d}}_q})\subseteq V_{\underline{d}}$ for each $\underline{d}\geq \underline{0}$ and $\widetilde{\underline{d}}_q\geq \underline{0}$.

Given a limit linear series, if $Y$ is an irreducible component of $X$, we denote by $V_{Y}$ the corresponding subspace of $H^{0}(X,\mathcal{L}_{Y})$. Also, we denote by $V_{il}$ the corresponding subspace of $H^{0}(X,\mathcal{L}_{il})$.

The conditions $\varphi_{\underline{d},\widetilde{\underline{d}}_q}(V_{\underline{d}})\subseteq V_{\widetilde{\underline{d}}_q}$ and $\varphi_{\widetilde{\underline{d}}_q,\underline{d}}(V_{\widetilde{\underline{d}}_q})\subseteq V_{\underline{d}}$ are called {\em the linking condition}, and we say that $V_{\underline{d}}$ and $V_{\widetilde{\underline{d}}_q}$ are linked by the maps $\varphi_{\underline{d},\widetilde{\underline{d}}_q}$ and $\varphi_{\widetilde{\underline{d}}_q,\underline{d}}$. 

Note that $\varphi_{\underline{d},\widetilde{\underline{d}}_q}: V_{\underline{d}}\rightarrow V_{\widetilde{\underline{d}}_q}$ has kernel $V_{\underline{d}}^{X_{q}^{c},0}$ and image contained in $V_{\widetilde{\underline{d}}_q}^{X_{q},0}$. Analogously, the map $\varphi_{\widetilde{\underline{d}}_q,\underline{d}}: V_{\widetilde{\underline{d}}_q}\rightarrow V_{\underline{d}}$ has kernel $V_{\widetilde{\underline{d}}_q}^{X_{q},0}$ and image contained in $V_{\underline{d}}^{X_{q}^{c},0}$.

A limit linear series $\{(\mathcal{L}_{\underline{d}},V_{\underline{d}})\}_{\underline{d}}$ is called {\em exact} if
\begin{center}
Im\,$(\varphi_{\underline{d},\widetilde{\underline{d}}_q}: V_{\underline{d}}\rightarrow V_{\widetilde{\underline{d}}_q})=V_{\widetilde{\underline{d}}_q}^{X_{q},0}$ and 
Im\,$(\varphi_{\widetilde{\underline{d}}_q,\underline{d}}: V_{\widetilde{\underline{d}}_q}\rightarrow V_{\underline{d}})=V_{\underline{d}}^{X_{q}^{c},0}$
\end{center}
for each $\underline{d}:=(i,d-i-l,l)\geq \underline{0}$ and $\widetilde{\underline{d}}_q\geq \underline{0}$.

}
\end{subs}

\begin{rem}\label{anulamiento}
{\em We have that, for $\underline{d}:=(i,d-i-l,l)\geq \underline{0}$, $\widetilde{\underline{d}}_q\geq \underline{0}$ and $\widetilde{q}\neq q$, $s\in H^{0}(\mathcal{L}_{\underline{d}})$ vanishes on $X_{\widetilde{q}}$ if and only if $\varphi_{\underline{d},\widetilde{\underline{d}}_q}(s)\in H^{0}(\mathcal{L}_{\widetilde{\underline{d}}_q})$ vanishes on $X_{\widetilde{q}}$. 
}
\end{rem}
Let us visualize the nonnegative multidegrees of total degree $d$ in the following ordered way:

\[\begin{array}{rcccc} (d,0,0) & (d-1,1,0) & \ldots & (1,d-1,0) & (0,d,0)\\
 & (d-1,0,1) & \ldots & (1,d-2,1) & (0,d-1,1)\\
 & & \ldots & \ldots & \ldots \\
 & &  & (1,0,d-1) & (0,1,d-1) \\
 & & & & (0,0,d)
\end{array}\]

Notice that multidegrees with the same first component are in the same column, multidegrees with the same third component are in the same row, the first component of the multidegrees increases from right to left, and the third component increases from top to bottom. We can use this order to visualize the relationship between different multidegrees. For instance, we have the following diagram:
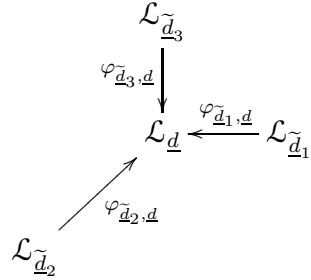
\[\xymatrix{&\mathcal{L}_{\widetilde{\underline{d}}_{3}}\ar[d]_{\varphi_{\widetilde{\underline{d}}_{3},\underline{d}}}&\\
&\mathcal{L}_{\underline{d}}&\mathcal{L}_{\widetilde{\underline{d}}_{1}}\ar[l]_{\varphi_{\widetilde{\underline{d}}_{1},\underline{d}}}\\
\mathcal{L}_{\widetilde{\underline{d}}_{2}}\ar[ru]_{\varphi_{\widetilde{\underline{d}}_{2},\underline{d}}}& &}\]

\section{Simple limit linear series}\label{simple lls}

Let $\underline{d}^{0},\underline{d}^{1},\ldots,\underline{d}^{n-1},\underline{d}$ be a sequence of nonnegative multidegrees of total degree $d$, with $\underline{d}^{0}\neq \underline{d}$. Put $\underline{d}^{n}:=\underline{d}$. Suppose the sequence satisfies the following condition:\\
For each $1\leq m\leq n$, 
\begin{center}
$\underline{d}^{m}=\underline{d}^{m-1}+\text{multideg}(\mathcal{O}_{(X_{q})})$ or $\underline{d}^{m-1}=\underline{d}^{m}+\text{multideg}(\mathcal{O}_{(X_{q})})$ for some $q$.
\end{center}
We say that the sequence $\underline{d}^{0},\ldots,\underline{d}$ forms a {\em path} connecting $\underline{d}^{0}$ to $\underline{d}$.

Given such a path, let $\varphi=\varphi_{\underline{d}^{n-1},\underline{d}}\circ\ldots\circ \varphi_{\underline{d}^{0},\underline{d}^{1}}$. Notice that $\varphi=0$ if at least one of the following conditions holds:\\
(I) There are multidegrees $\underline{d}^{m_{1}},\underline{d}^{m_{2}}$ and $\underline{d}^{m_{3}}$, with $m_{q}<n$ for $q=1,2,3$, such that
\begin{center}
$\underline{d}^{m_{q}+1}=\underline{d}^{m_{q}}+\text{multideg}(\mathcal{O}_{(X_{q})})$ for each $q$.
\end{center}
(II) There are multidegrees $\underline{d}^{m_{1}}$ and $\underline{d}^{m_{2}}$, with $m_{p}<n$ for $p=1,2$, such that
\begin{center}
$\underline{d}^{m_{1}+1}=\underline{d}^{m_{1}}+\text{multideg}(\mathcal{O}_{(X_{q})})$ and $\underline{d}^{m_{2}}=\underline{d}^{m_{2}+1}+\text{multideg}(\mathcal{O}_{(X_{q})})$ for some $q$.
\end{center}
(III) There are multidegrees $\underline{d}^{m_{1}}$ and $\underline{d}^{m_{2}}$, with $m_{p}<n$ for $p=1,2$, such that
\begin{center}
$\underline{d}^{m_{1}}=\underline{d}^{m_{1}+1}+\text{multideg}(\mathcal{O}_{(X_{q_{1}})})$ and $\underline{d}^{m_{2}}=\underline{d}^{m_{2}+1}+\text{multideg}(\mathcal{O}_{(X_{q_{2}})})$ for some $q_{1}\neq q_{2}$.
\end{center}
\begin{rem}\label{path existence}
{\em Given two multidegrees $\underline{d}^{0}=(\tilde{i},d-\tilde{i}-\tilde{l},\tilde{l})\geq \underline{0}$ and $\underline{d}=(i,d-i-l,l)\geq \underline{0}$, there exist nonnegative multidegrees $\underline{d}^{1},\ldots,\underline{d}^{n-1}$ of total degree $d$ such that the sequence $\underline{d}^{0},\underline{d}^{1},\ldots,\underline{d}^{n-1},\underline{d}^{n}:=\underline{d}$ forms a path connecting $\underline{d}^{0}$ to $\underline{d}$ and does not satisfy any of the conditions (I), (II)and (III). The proof is straightforward. We will only see two cases, as the other cases are analogous. If $i<\tilde{i}$ and $l>\tilde{l}$, we can consider the horizontal path connecting $\underline{d}^{0}=(\tilde{i},d-\tilde{i}-\tilde{l},\tilde{l})$ to $(i,d-i-\tilde{l},\tilde{l})$ followed by the vertical path connecting $(i,d-i-\tilde{l},\tilde{l})$ to $\underline{d}=(i,d-i-l,l)$, that is
\begin{center}
$\underline{d}^{1}=(\tilde{i}-1,d-(\tilde{i}-1)-\tilde{l},\tilde{l}),\ldots,\underline{d}^{\tilde{i}-i}=(i,d-i-\tilde{l},\tilde{l})$,\\
$\underline{d}^{\tilde{i}-i+1}=(i,d-i-(\tilde{l}+1),\tilde{l}+1),\ldots,\underline{d}^{\tilde{i}-i+l-\tilde{l}}=(i,d-i-l,l)$.
\end{center}
If $\tilde{i}<i$, $\tilde{l}<l$ and $i-\tilde{i}<l-\tilde{l}$, we can consider the diagonal path connecting the multidegree $\underline{d}^{0}=(\tilde{i},d-\tilde{i}-\tilde{l},\tilde{l})$ to $(i,d-i-(\tilde{l}+i-\tilde{i}),\tilde{l}+i-\tilde{i})$ followed by the vertical path connecting $(i,d-i-(\tilde{l}+i-\tilde{i}),\tilde{l}+i-\tilde{i})$ to $\underline{d}=(i,d-i-l,l)$, that is
\begin{center}
$\underline{d}^{1}=(\tilde{i}+1,d-(\tilde{i}+1)-(\tilde{l}+1),\tilde{l}+1),\ldots,\underline{d}^{i-\tilde{i}}=(i,d-i-(\tilde{l}+i-\tilde{i}),\tilde{l}+i-\tilde{i})$,\\
$\underline{d}^{i-\tilde{i}+1}=(i,d-i-(\tilde{l}+i-\tilde{i}+1),\tilde{l}+i-\tilde{i}+1),\ldots,\underline{d}^{l-\tilde{l}}=(i,d-i-l,l)$.
\end{center}
}
\end{rem}

\begin{rem}\label{path independence}
{\em Let $\underline{d}^{0}=(\tilde{i},d-\tilde{i}-\tilde{l},\tilde{l})\geq \underline{0}$ and $\underline{d}=(i,d-i-l,l)\geq \underline{0}$, with $\underline{d}^{0}\neq \underline{d}$. Then, the maps $\varphi_{\underline{d}^{n-1},\underline{d}}\circ\ldots\circ \varphi_{\underline{d}^{0},\underline{d}^{1}}$ obtained from sequences of nonnegative multidegrees $\underline{d}^{0},\underline{d}^{1},\ldots,\underline{d}^{n-1},\underline{d}^{n}:=\underline{d}$ of total degree $d$ forming a path connecting $\underline{d}^{0}$ to $\underline{d}$ and not satisfying any of the conditions (I), (II) and (III) are equal because of the property (3) of the maps $\varphi_{\underline{d},\widetilde{\underline{d}}_q}$ and $\varphi_{\widetilde{\underline{d}}_q,\underline{d}}$ we have seen in Subsection \ref{limit1}. We put $\varphi_{\underline{d}^{0},\underline{d}}=\varphi_{\underline{d}^{n-1},\underline{d}}\circ\ldots\circ \varphi_{\underline{d}^{0},\underline{d}^{1}}$ for any such a sequence. Also, we put $\varphi_{\underline{d},\underline{d}}=\text{id}$ for any $\underline{d}$.
}
\end{rem}

\begin{defn}\label{simple lls definition}
{\em Let $\{(\mathcal{L}_{\underline{d}},V_{\underline{d}})\}_{\underline{d}}$ be a limit linear series of degree $d$ and dimension $r$. We say that $\{(\mathcal{L}_{\underline{d}},V_{\underline{d}})\}_{\underline{d}}$ is {\em simple} if there exist $m$ nonnegative multidegrees $\underline{d}^{1},\ldots,\underline{d}^{m}$ of total degree $d$ and linearly independent sections $s^{\alpha}_{1},\ldots,s^{\alpha}_{r_{\alpha}}\in V_{\underline{d}^{\alpha}}$ for each $\alpha=1,\ldots,m$ such that $\varphi_{\underline{d}^{1},\underline{d}}(s^{1}_{1}),\ldots,\varphi_{\underline{d}^{1},\underline{d}}(s^{1}_{r_{1}}),\ldots,\varphi_{\underline{d}^{m},\underline{d}}(s^{m}_{1}),\ldots,\varphi_{\underline{d}^{m},\underline{d}}(s^{m}_{r_{m}})$ form a basis for $V_{\underline{d}}$, for each $\underline{d}$.
}
\end{defn}

\begin{prop}\label{simple is exact}
Let $\{(\mathcal{L}_{\underline{d}},V_{\underline{d}})\}_{\underline{d}}$ be a simple limit linear series. Then it is exact.
\end{prop}
{\em Proof.} We will first see how the horizontal exactness and the vertical exactness imply the diagonal exactness. In fact, for $i\geq 1$ and $l\geq 1$ such that $i+l\leq d$, let $\underline{d}=(i,d-i-l,l)$, $\underline{d}'=(i-1,d-i-l+2,l-1)$ and $\underline{d}''=(i-1,d-i-l+1,l)$. Suppose $\varphi_{\underline{d},\underline{d}''}(V_{\underline{d}})=V_{\underline{d}''}^{X_{1},0}$ and $\varphi_{\underline{d}'',\underline{d}'}(V_{\underline{d}''})=V_{\underline{d}'}^{X_{3},0}$. Then
\begin{align}
\varphi_{\underline{d},\underline{d}'}(V_{\underline{d}})=&(\varphi_{\underline{d}'',\underline{d}'}\circ \varphi_{\underline{d},\underline{d}''})(V_{\underline{d}}) \nonumber \\
=&\varphi_{\underline{d}'',\underline{d}'}(\varphi_{\underline{d},\underline{d}''}(V_{\underline{d}})) \nonumber \\
=&\varphi_{\underline{d}'',\underline{d}'}(V_{\underline{d}''}^{X_{1},0}) \nonumber \\
=&(\varphi_{\underline{d}'',\underline{d}'}(V_{\underline{d}''}))^{X_{1},0}, \nonumber
\end{align}
where in the last equality we used Remark \ref{anulamiento}. Thus
\begin{equation}\label{simple is exact1}
\varphi_{\underline{d},\underline{d}'}(V_{\underline{d}})=(\varphi_{\underline{d}'',\underline{d}'}(V_{\underline{d}''}))^{X_{1},0}=(V_{\underline{d}'}^{X_3,0})^{X_{1},0}=V_{\underline{d}'}^{X_{2}^{c},0}.
\end{equation}
On the other hand,
\begin{center}
$\text{dim}\,V_{\underline{d}'}^{X_{2}^{c},0}+\text{dim}\,\varphi_{\underline{d}',\underline{d}}(V_{\underline{d}'})=\text{dim}\,V_{\underline{d}'}=r+1$,
\end{center}
so
\begin{equation*}
\text{dim}\,\varphi_{\underline{d}',\underline{d}}(V_{\underline{d}'})=r+1-\text{dim}\,V_{\underline{d}'}^{X_{2}^{c},0}
=\text{dim}\,V_{\underline{d}}-\text{dim}\,V_{\underline{d}'}^{X_{2}^{c},0}
=\text{dim}\,V_{\underline{d}}^{X_{2},0},
\end{equation*}
where the last equality follows from (\ref{simple is exact1}). Since $\varphi_{\underline{d}',\underline{d}}(V_{\underline{d}'})\subseteq V_{\underline{d}}^{X_{2},0}$, it follows that 
\begin{equation}\label{simple is exact2}
\varphi_{\underline{d}',\underline{d}}(V_{\underline{d}'})=V_{\underline{d}}^{X_{2},0}
\end{equation}
Thus, the diagonal exactness follows from (\ref{simple is exact1}) and (\ref{simple is exact2}).

It remains to show the horizontal and vertical exactness. We will only prove the horizontal exactness, as the proof of the vertical exactness is analogous. For $i\geq 1$ and $l\geq 0$ such that $i+l\leq d$, let $\underline{d}=(i,d-i-l,l)$. Then $\widetilde{\underline{d}}_{1}=(i-1,d-i-l+1,l)\geq \underline{0}$.

Let $M$ be the set of nonnegative multidegrees of total degree $d$. We define the following subsets of $M$:
\begin{center}
$D=\{\underline{d}_{0}=(\tilde{i},d-\tilde{i}-\tilde{l},\tilde{l})\in M\,/\,\,\tilde{l}\leq l\,\,\text{and}\,\,\tilde{i}-i\geq \tilde{l}-l\}$,\\
$E=\{\underline{d}_{0}=(\tilde{i},d-\tilde{i}-\tilde{l},\tilde{l})\in M\,/\,\,\tilde{l}\geq l\,\,\text{and}\,\,\tilde{i}\geq i\}$ and\\
$F=\{\underline{d}_{0}=(\tilde{i},d-\tilde{i}-\tilde{l},\tilde{l})\in M\,/\,\,\tilde{i}\leq i-1\,\,\text{and}\,\,\tilde{i}-(i-1)\leq \tilde{l}-l\}$
\end{center}
Notice that $(D\cup E)\cap F=\emptyset$ and $M=(D\cup E)\cup F$.

Since $\{(\mathcal{L}_{\underline{d}},V_{\underline{d}})\}_{\underline{d}}$ is simple, there exist multidegrees $\underline{d}^{1},\ldots,\underline{d}^{m}$ and linearly independent sections $s^{\alpha}_{1},\ldots,s^{\alpha}_{r_{\alpha}}\in V_{\underline{d}^{\alpha}}$ for each $\alpha$ as in Definition \ref{simple lls definition}. We claim that 
\begin{equation}\label{composition1}
\varphi_{\widetilde{\underline{d}}_{1},\underline{d}}\circ \varphi_{\underline{d}^{\alpha},\widetilde{\underline{d}}_{1}}=\varphi_{\underline{d}^{\alpha},\underline{d}},
\end{equation}
for any $\underline{d}^{\alpha}\in F$. In fact, let $\underline{d}^{\alpha}=(\tilde{i},d-\tilde{i}-\tilde{l},\tilde{l})\in F$. If $\tilde{l}\geq l$, we can obtain $\varphi_{\underline{d}^{\alpha},\widetilde{\underline{d}}_{1}}$ by considering the vertical path connecting $\underline{d}^{\alpha}=(\tilde{i},d-\tilde{i}-\tilde{l},\tilde{l})$ to $(\tilde{i},d-\tilde{i}-l,l)$ followed by the horizontal path connecting $(\tilde{i},d-\tilde{i}-l,l)$ to $\widetilde{\underline{d}}_{1}=(i-1,d-i-l+1,l)$. Analogously, we can obtain $\varphi_{\underline{d}^{\alpha},\underline{d}}$ by considering the vertical path connecting $\underline{d}^{\alpha}=(\tilde{i},d-\tilde{i}-\tilde{l},\tilde{l})$ to $(\tilde{i},d-\tilde{i}-l,l)$ followed by the horizontal path connecting $(\tilde{i},d-\tilde{i}-l,l)$ to $\underline{d}=(i,d-i-l,l)$. Since $\tilde{i}\leq i-1$, we conclude that
\begin{center}
$\varphi_{\underline{d}^{\alpha},\underline{d}}=\varphi_{\widetilde{\underline{d}}_{1},\underline{d}}\circ \varphi_{\underline{d}^{\alpha},\widetilde{\underline{d}}_{1}}$.
\end{center}
If $\tilde{l}<l$, we can obtain $\varphi_{\underline{d}^{\alpha},\widetilde{\underline{d}}_{1}}$ by considering the diagonal path connecting the multidegree $\underline{d}^{\alpha}=(\tilde{i},d-\tilde{i}-\tilde{l},\tilde{l})$ to $(\tilde{i}+l-\tilde{l},d-(\tilde{i}+l-\tilde{l})-l,l)$ followed by the horizontal path connecting $(\tilde{i}+l-\tilde{l},d-(\tilde{i}+l-\tilde{l})-l,l)$ to $\widetilde{\underline{d}}_{1}=(i-1,d-i-l+1,l)$. Analogously, we can obtain $\varphi_{\underline{d}^{\alpha},\underline{d}}$ by considering the diagonal path connecting $\underline{d}^{\alpha}=(\tilde{i},d-\tilde{i}-\tilde{l},\tilde{l})$ to $(\tilde{i}+l-\tilde{l},d-(\tilde{i}+l-\tilde{l})-l,l)$ followed by the horizontal path connecting $(\tilde{i}+l-\tilde{l},d-(\tilde{i}+l-\tilde{l})-l,l)$ to $\underline{d}=(i,d-i-l,l)$. Since $\tilde{i}+l-\tilde{l}\leq i-1$, we conclude that
\begin{center}
$\varphi_{\underline{d}^{\alpha},\underline{d}}=\varphi_{\widetilde{\underline{d}}_{1},\underline{d}}\circ \varphi_{\underline{d}^{\alpha},\widetilde{\underline{d}}_{1}}$.
\end{center}
This proves the claim (\ref{composition1}).

Now, we claim that 
\begin{equation}\label{composition2}
\varphi_{\underline{d},\widetilde{\underline{d}}_{1}}\circ \varphi_{\underline{d}^{\alpha},\underline{d}}=\varphi_{\underline{d}^{\alpha},\widetilde{\underline{d}}_{1}},
\end{equation}
for any $\underline{d}^{\alpha}\in D\cup E$. In fact, let $\underline{d}^{\alpha}=(\tilde{i},d-\tilde{i}-\tilde{l},\tilde{l})\in D\cup E$. If $\underline{d}^{\alpha}\in E$, the proof of (\ref{composition2}) is analogous to that of the claim (\ref{composition1}) for the case $\tilde{l}\geq l$.

Suppose now that $\underline{d}^{\alpha}\in D$. If $\tilde{i}\geq i$, again, the proof of (\ref{composition2}) is analogous to that of the claim (\ref{composition1}) for the first case.

If $\tilde{i}<i$, we can obtain $\varphi_{\underline{d}^{\alpha},\underline{d}}$ by considering the diagonal path connecting the multidegree $\underline{d}^{\alpha}=(\tilde{i},d-\tilde{i}-\tilde{l},\tilde{l})$ to $(i,d-i-(\tilde{l}+i-\tilde{i}),\tilde{l}+i-\tilde{i})$ followed by the vertical path connecting $(i,d-i-(\tilde{l}+i-\tilde{i}),\tilde{l}+i-\tilde{i})$ to $\underline{d}=(i,d-i-l,l)$. Call $\gamma$ this path considered to obtain $\varphi_{\underline{d}^{\alpha},\underline{d}}$. Notice that $\gamma$ followed by the arrow from $\underline{d}$ to $\widetilde{\underline{d}}_{1}$ is a path, which does not satisfy any of the conditions (I), (II) and (III). We conclude that  
\begin{center}
$\varphi_{\underline{d},\widetilde{\underline{d}}_{1}}\circ \varphi_{\underline{d}^{\alpha},\underline{d}}=\varphi_{\underline{d}^{\alpha},\widetilde{\underline{d}}_{1}}$.
\end{center}
This proves the claim (\ref{composition2}).

On the other hand, we have that 
\begin{center}
$B_{\underline{d}}:=\{\varphi_{\underline{d}^{1},\underline{d}}(s^{1}_{1}),\ldots,\varphi_{\underline{d}^{1},\underline{d}}(s^{1}_{r_{1}}),\ldots,\varphi_{\underline{d}^{m},\underline{d}}(s^{m}_{1}),\ldots,\varphi_{\underline{d}^{m},\underline{d}}(s^{m}_{r_{m}})\}$
\end{center}
is a basis of $V_{\underline{d}}$. Now, by the claim (\ref{composition1}), for each $\underline{d}^{\alpha}\in F$ we have 
\begin{align}
\varphi_{\underline{d},\widetilde{\underline{d}}_{1}}\circ \varphi_{\underline{d}^{\alpha},\underline{d}}=&\varphi_{\underline{d},\widetilde{\underline{d}}_{1}}\circ (\varphi_{\widetilde{\underline{d}}_{1},\underline{d}}\circ \varphi_{\underline{d}^{\alpha},\widetilde{\underline{d}}_{1}})  \nonumber \\
=&(\varphi_{\underline{d},\widetilde{\underline{d}}_{1}}\circ \varphi_{\widetilde{\underline{d}}_{1},\underline{d}})\circ \varphi_{\underline{d}^{\alpha},\widetilde{\underline{d}}_{1}} \nonumber \\
=&0, \nonumber
\end{align}
as $\varphi_{\underline{d},\widetilde{\underline{d}}_{1}}\circ \varphi_{\widetilde{\underline{d}}_{1},\underline{d}}=0$. Then, we have 
\begin{equation}\label{inclusion}
B_{\underline{d},F}:=\{\varphi_{\underline{d}^{\alpha},\underline{d}}(s^{\alpha}_{z})\,/\,\,\underline{d}^{\alpha}\in F\,\text{and}\,z=1,\ldots,r_{\alpha}\}\subseteq \text{Ker}\,(\varphi_{\underline{d},\widetilde{\underline{d}}_{1}}:V_{\underline{d}}\rightarrow V_{\widetilde{\underline{d}}_{1}}).
\end{equation}
On the other hand, we have that 
\begin{center}
$B_{\widetilde{\underline{d}}_{1}}:=\{\varphi_{\underline{d}^{1},\widetilde{\underline{d}}_{1}}(s^{1}_{1}),\ldots,\varphi_{\underline{d}^{1},\widetilde{\underline{d}}_{1}}(s^{1}_{r_{1}}),\ldots,\varphi_{\underline{d}^{m},\widetilde{\underline{d}}_{1}}(s^{m}_{1}),\ldots,\varphi_{\underline{d}^{m},\widetilde{\underline{d}}_{1}}(s^{m}_{r_{m}})\}$
\end{center}
is a basis of $V_{\widetilde{\underline{d}}_{1}}$. In particular,
\begin{center}
$\{\varphi_{\underline{d}^{\alpha},\widetilde{\underline{d}}_{1}}(s^{\alpha}_{z})\,/\,\,\underline{d}^{\alpha}\in D\cup E\,\text{and}\,z=1,\ldots,r_{\alpha}\}$ is linearly independent.
\end{center}
By the claim (\ref{composition2}), this means that
\begin{center}
$\{\varphi_{\underline{d},\widetilde{\underline{d}}_{1}}(\varphi_{\underline{d}^{\alpha},\underline{d}}(s^{\alpha}_{z}))\,/\,\,\underline{d}^{\alpha}\in D\cup E\,\text{and}\,z=1,\ldots,r_{\alpha}\}$ is linearly independent,
\end{center}
and hence
\begin{equation}\label{intersection}
\left\langle B_{\underline{d},D\cup E}\right\rangle\cap \text{Ker}\,(\varphi_{\underline{d},\widetilde{\underline{d}}_{1}}:V_{\underline{d}}\rightarrow V_{\widetilde{\underline{d}}_{1}})=0,
\end{equation}
where $B_{\underline{d},D\cup E}:=\{\varphi_{\underline{d}^{\alpha},\underline{d}}(s^{\alpha}_{z})\,/\,\,\underline{d}^{\alpha}\in D\cup E\,\text{and}\,z=1,\ldots,r_{\alpha}\}$. Also, notice that the basis $B_{\underline{d}}$ de $V_{\underline{d}}$ can be written as
\begin{equation}\label{basis union}
B_{\underline{d}}=B_{\underline{d},F}\cup B_{\underline{d},D\cup E},
\end{equation}
as $M=(D\cup E)\cup F$. It follows from (\ref{inclusion}), (\ref{intersection}) and (\ref{basis union}) that
\begin{center}
$\text{Ker}\,(\varphi_{\underline{d},\widetilde{\underline{d}}_{1}}:V_{\underline{d}}\rightarrow V_{\widetilde{\underline{d}}_{1}})=\left\langle B_{\underline{d},F}\right\rangle$.
\end{center}
But, it follows from the claim (\ref{composition1}) that $B_{\underline{d},F}\subseteq \text{Im}\,(\varphi_{\widetilde{\underline{d}}_{1},\underline{d}}:V_{\widetilde{\underline{d}}_{1}}\rightarrow V_{\underline{d}})$. Thus
\begin{center}
$\text{Ker}\,(\varphi_{\underline{d},\widetilde{\underline{d}}_{1}}:V_{\underline{d}}\rightarrow V_{\widetilde{\underline{d}}_{1}})\subseteq \text{Im}\,(\varphi_{\widetilde{\underline{d}}_{1},\underline{d}}:V_{\widetilde{\underline{d}}_{1}}\rightarrow V_{\underline{d}})$,
\end{center}
and hence
\begin{center}
$\text{Ker}\,(\varphi_{\underline{d},\widetilde{\underline{d}}_{1}}:V_{\underline{d}}\rightarrow V_{\widetilde{\underline{d}}_{1}})=\text{Im}\,(\varphi_{\widetilde{\underline{d}}_{1},\underline{d}}:V_{\widetilde{\underline{d}}_{1}}\rightarrow V_{\underline{d}})$.
\end{center}
The proof of 
\begin{center}
$\text{Ker}\,(\varphi_{\widetilde{\underline{d}}_{1},\underline{d}}:V_{\widetilde{\underline{d}}_{1}}\rightarrow V_{\underline{d}})=\text{Im}\,(\varphi_{\underline{d},\widetilde{\underline{d}}_{1}}:V_{\underline{d}}\rightarrow V_{\widetilde{\underline{d}}_{1}})$
\end{center}
is analogous. This proves the horizontal exactness, which finishes the proof of the proposition. 
\hfill $\Box$

Now, we will show some results which will be useful to obtain some properties of simple limit linear series.
\begin{prop}\label{codimension1}
For any $i\geq 0$ and $l\geq 0$ such that $i+l\leq d$, let $\underline{d}:=(i,d-i-l,l)$. Then, the following statements hold:\\
1. If $i>0$ and $l>0$, let $\underline{d}':=(i-1,d-i-l+2,l-1)$. Let $V_{\underline{d}}$ and $V_{\underline{d}'}$ be $r+1$-dimensional subspaces of $H^{0}(\mathcal{L}_{\underline{d}})$ and $H^{0}(\mathcal{L}_{\underline{d}'})$, respectively, such that $\varphi_{\underline{d}',\underline{d}}(V_{\underline{d}'})=V_{\underline{d}}^{X_{2},0}$. Then
\begin{center}
$r+1\,-$\,{\em dim}\,$(V_{\underline{d}'}^{X_{1},0}+V_{\underline{d}'}^{X_{3},0})=${\em dim}\,$V_{\underline{d}}^{X_{2},0}\,-$\,{\em dim}\,$(V_{\underline{d}}^{X_{3}^{c},0}\oplus V_{\underline{d}}^{X_{1}^{c},0})$.
\end{center}
2. If $i>0$, let $\underline{d}'':=(i-1,d-i-l+1,l)$. Let $V_{\underline{d}}$ and $V_{\underline{d}''}$ be $r+1$-dimensional subspaces of $H^{0}(\mathcal{L}_{\underline{d}})$ and $H^{0}(\mathcal{L}_{\underline{d}''})$, respectively, such that $\varphi_{\underline{d},\underline{d}''}(V_{\underline{d}})=V_{\underline{d}''}^{X_{1},0}$. Then 
\begin{center}
$r+1\,-$\,{\em dim}\,$(V_{\underline{d}}^{X_{2},0}+V_{\underline{d}}^{X_{3},0})=${\em dim}\,$V_{\underline{d}''}^{X_{1},0}\,-$\,{\em dim}\,$(V_{\underline{d}''}^{X_{3}^{c},0}\oplus V_{\underline{d}''}^{X_{2}^{c},0})$.
\end{center}
3. If $i+l\leq d-1$, let $\underline{d}''':=(i,d-i-l-1,l+1)$. Let $V_{\underline{d}}$ and $V_{\underline{d}'''}$ be $r+1$-dimensional subspaces of $H^{0}(\mathcal{L}_{\underline{d}})$ and $H^{0}(\mathcal{L}_{\underline{d}'''})$, respectively, such that $\varphi_{\underline{d}''',\underline{d}}(V_{\underline{d}'''})=V_{\underline{d}}^{X_{3},0}$. Then
\begin{center}
$r+1\,-$\,{\em dim}\,$(V_{\underline{d}'''}^{X_{1},0}+V_{\underline{d}'''}^{X_{2},0})=${\em dim}\,$V_{\underline{d}}^{X_{3},0}\,-$\,{\em dim}\,$(V_{\underline{d}}^{X_{2}^{c},0}\oplus V_{\underline{d}}^{X_{1}^{c},0})$.
\end{center}
\end{prop}
{\em Proof.} We will only prove the statement 2, as the proofs of the other statements are analogous. Set $\beta:=r+1-\text{dim}\,(V_{\underline{d}}^{X_{2},0}+V_{\underline{d}}^{X_{3},0})$ and write
\begin{center}
$V_{\underline{d}}=(V_{\underline{d}}^{X_{2},0}+V_{\underline{d}}^{X_{3},0})\oplus \left\langle v_1,\ldots,v_{\beta}\right\rangle$,
\end{center}
for some $v_1,\ldots,v_{\beta}$ which are linearly independent. Since
\begin{center}
$V_{\underline{d}}^{X_{2},0}+V_{\underline{d}}^{X_{3},0}\supseteq V_{\underline{d}}^{X_{1}^{c},0}=\text{Ker}\,(\varphi_{\underline{d},\underline{d}''}|_{V_{\underline{d}}})$,
\end{center}
we conclude that $\varphi_{\underline{d},\underline{d}''}(v_1),\ldots,\varphi_{\underline{d},\underline{d}''}(v_{\beta})$ are linearly independent and
\begin{center}
$\varphi_{\underline{d},\underline{d}''}(V_{\underline{d}})=(\varphi_{\underline{d},\underline{d}''}(V_{\underline{d}}^{X_{2},0})+\varphi_{\underline{d},\underline{d}''}(V_{\underline{d}}^{X_{3},0}))\oplus \left\langle \varphi_{\underline{d},\underline{d}''}(v_1),\ldots,\varphi_{\underline{d},\underline{d}''}(v_{\beta})\right\rangle$.
\end{center}
On the other hand, as $\varphi_{\underline{d},\underline{d}''}(V_{\underline{d}})=V_{\underline{d}''}^{X_{1},0}$, it follows that
$\varphi_{\underline{d},\underline{d}''}(V_{\underline{d}}^{X_{2},0})=V_{\underline{d}''}^{X_{3}^{c},0}$ and $\varphi_{\underline{d},\underline{d}''}(V_{\underline{d}}^{X_{3},0})=V_{\underline{d}''}^{X_{2}^{c},0}$. Thus 
\begin{center}
$V_{\underline{d}''}^{X_{1},0}=(V_{\underline{d}''}^{X_{3}^{c},0}+V_{\underline{d}''}^{X_{2}^{c},0})\oplus \left\langle \varphi_{\underline{d},\underline{d}''}(v_1),\ldots,\varphi_{\underline{d},\underline{d}''}(v_{\beta})\right\rangle$.
\end{center}
As $\varphi_{\underline{d},\underline{d}''}(v_1),\ldots,\varphi_{\underline{d},\underline{d}''}(v_{\beta})$ are linearly independent, it follows that
\begin{center}
$\beta=\text{dim}\,V_{\underline{d}''}^{X_{1},0}-\text{dim}\,(V_{\underline{d}''}^{X_{3}^{c},0}+V_{\underline{d}''}^{X_{2}^{c},0})$,
\end{center}
which proves the statement 2.
\hfill $\Box$

From now on we will keep the notation of multidegrees used in Proposition \ref{codimension1}.

\begin{prop}\label{suma1}
Let $V_{\underline{d}}$ and $V_{\underline{d}'''}$ be $r+1$-dimensional subspaces of $H^{0}(\mathcal{L}_{\underline{d}})$ and $H^{0}(\mathcal{L}_{\underline{d}'''})$, respectively, such that $\varphi_{\underline{d}''',\underline{d}}(V_{\underline{d}'''})=V_{\underline{d}}^{X_{3},0}$. If $V_{\underline{d}'''}=V_{\underline{d}'''}^{X_{2},0}\oplus \left\langle v_1,\ldots,v_{m}\right\rangle$, where $v_1,\ldots,v_{m}$ are linearly independent, then $\varphi_{\underline{d}''',\underline{d}}(v_1),\ldots,\varphi_{\underline{d}''',\underline{d}}(v_{m})$ are linearly independent and
\begin{center}
$V_{\underline{d}}^{X_{2},0}+V_{\underline{d}}^{X_{3},0}=V_{\underline{d}}^{X_{2},0}\oplus \left\langle \varphi_{\underline{d}''',\underline{d}}(v_1),\ldots,\varphi_{\underline{d}''',\underline{d}}(v_{m})\right\rangle$.
\end{center}
\end{prop}
{\em Proof.} Since 
\begin{center}
$V_{\underline{d}'''}^{X_{2},0}\supseteq V_{\underline{d}'''}^{X_{3}^{c},0}=\text{Ker}\,(\varphi_{\underline{d}''',\underline{d}}|_{V_{\underline{d}'''}})$ and $V_{\underline{d}'''}=V_{\underline{d}'''}^{X_{2},0}\oplus \left\langle v_1,\ldots,v_{m}\right\rangle$,
\end{center}
it follows that $\varphi_{\underline{d}''',\underline{d}}(v_1),\ldots,\varphi_{\underline{d}''',\underline{d}}(v_{m})$ are linearly independent and
\begin{center}
$\varphi_{\underline{d}''',\underline{d}}(V_{\underline{d}'''})=\varphi_{\underline{d}''',\underline{d}}(V_{\underline{d}'''}^{X_{2},0})\oplus \left\langle \varphi_{\underline{d}''',\underline{d}}(v_1),\ldots,\varphi_{\underline{d}''',\underline{d}}(v_{m})\right\rangle$.
\end{center}
On the other hand, as $\varphi_{\underline{d}''',\underline{d}}(V_{\underline{d}'''})=V_{\underline{d}}^{X_{3},0}$, it follows that $\varphi_{\underline{d}''',\underline{d}}(V_{\underline{d}'''}^{X_{2},0})=V_{\underline{d}}^{X_{1}^{c},0}$. Thus
\begin{center}
$V_{\underline{d}}^{X_{3},0}=V_{\underline{d}}^{X_{1}^{c},0}\oplus \left\langle \varphi_{\underline{d}''',\underline{d}}(v_1),\ldots,\varphi_{\underline{d}''',\underline{d}}(v_{m})\right\rangle$.
\end{center}
Now, since $V_{\underline{d}}^{X_{1}^{c},0}=V_{\underline{d}}^{X_{2},0}\cap V_{\underline{d}}^{X_{3},0}$, we conclude that
\begin{center}
$V_{\underline{d}}^{X_{2},0}+V_{\underline{d}}^{X_{3},0}=V_{\underline{d}}^{X_{2},0}\oplus \left\langle \varphi_{\underline{d}''',\underline{d}}(v_1),\ldots,\varphi_{\underline{d}''',\underline{d}}(v_{m})\right\rangle$.
\end{center}
\hfill $\Box$

\begin{cor}\label{suma1corolario1}
Let $V_{\underline{d}}$ and $V_{\underline{d}'''}$ be $r+1$-dimensional subspaces of $H^{0}(\mathcal{L}_{\underline{d}})$ and $H^{0}(\mathcal{L}_{\underline{d}'''})$, respectively, such that $\varphi_{\underline{d}''',\underline{d}}(V_{\underline{d}'''})=V_{\underline{d}}^{X_{3},0}$. Then
\begin{center}
{\em dim}\,$V_{\underline{d}'''}^{X_{2},0}\,-$\,{\em dim}\,$V_{\underline{d}}^{X_{2},0}=r+1\,-$\,{\em dim}\,$(V_{\underline{d}}^{X_{2},0}+V_{\underline{d}}^{X_{3},0})$.
\end{center}
\end{cor}
{\em Proof.} Write 
\begin{equation}\label{primera suma}
V_{\underline{d}'''}=V_{\underline{d}'''}^{X_{2},0}\oplus \left\langle v_1,\ldots,v_{m}\right\rangle, 
\end{equation}
where $v_1,\ldots,v_{m}$ are linearly independent. By Proposition \ref{suma1}, we have that 
\begin{center}
$\varphi_{\underline{d}''',\underline{d}}(v_1),\ldots,\varphi_{\underline{d}''',\underline{d}}(v_{m})$ are linearly independent and
\end{center}
\begin{equation}\label{segunda suma}
V_{\underline{d}}^{X_{2},0}+V_{\underline{d}}^{X_{3},0}=V_{\underline{d}}^{X_{2},0}\oplus \left\langle \varphi_{\underline{d}''',\underline{d}}(v_1),\ldots,\varphi_{\underline{d}''',\underline{d}}(v_{m})\right\rangle.
\end{equation}
It follows from (\ref{primera suma}) that
\begin{equation}\label{m despejado}
m=r+1-\text{dim}\,V_{\underline{d}'''}^{X_{2},0}
\end{equation}
Thus, from (\ref{segunda suma}) and (\ref{m despejado}), we get
\begin{align}
\text{dim}\,(V_{\underline{d}}^{X_{2},0}+V_{\underline{d}}^{X_{3},0})&=\text{dim}\,V_{\underline{d}}^{X_{2},0}+m \nonumber \\
&=\text{dim}\,V_{\underline{d}}^{X_{2},0}+r+1-\text{dim}\,V_{\underline{d}'''}^{X_{2},0}, \nonumber 
\end{align}
and hence
\begin{center}
$\text{dim}\,V_{\underline{d}'''}^{X_{2},0}-\text{dim}\,V_{\underline{d}}^{X_{2},0}=r+1-\text{dim}\,(V_{\underline{d}}^{X_{2},0}+V_{\underline{d}}^{X_{3},0})$.
\end{center}
\hfill $\Box$

\begin{rem}\label{lado diagonal}
{\em Given $i\geq 0$ and $l\geq 0$ such that $i+l=d$, and a subspace $V\subseteq H^{0}(\mathcal{L}_{il})$, we have that $V^{X_{1},0}\subseteq V^{X_{2},0}$ and $V^{X_{3},0}\subseteq V^{X_{2},0}$. We will prove the first inclusion, as the proof of the other inclusion is analogous. Let $s\in V^{X_{1},0}$. Then $s|_{X_{2}}\in H^{0}(\mathcal{L}_{il}|_{X_{2}}(-A))$. On the other hand, since $i+l=d$, we have $\text{multideg}(\mathcal{L}_{il})=(i,0,l)$, which implies that $\text{deg}(\mathcal{L}_{il}|_{X_{2}}(-A))=-1<0$. Then $H^{0}(\mathcal{L}_{il}|_{X_{2}}(-A))=0$, so $s|_{X_{2}}=0$. Thus $s\in V^{X_{2},0}$, which proves the desired inclusion.  
}
\end{rem}

\begin{rem}\label{vertice superior derecho}
{\em For any subspace $V\subseteq H^{0}(\mathcal{L}_{00})$, we have that $V^{X_{2},0}=0$. In fact, let $s\in V^{X_{2},0}$. Then $s|_{X_{1}}\in H^{0}(\mathcal{L}_{00}|_{X_{1}}(-A))$. On the other hand, since the multidegree of $\mathcal{L}_{00}$ is $(0,d,0)$, we have that $\text{deg}(\mathcal{L}_{00}|_{X_{1}}(-A))=-1<0$. Then $H^{0}(\mathcal{L}_{00}|_{X_{1}}(-A))=0$, so $s|_{X_{1}}=0$. Analogously, $s|_{X_{3}}=0$. Therefore $s=0$. We conclude that $V^{X_{2},0}=0$.
}
\end{rem}

\begin{prop}\label{dimension de cocientes}
The following statements hold:\\
1. Let $V_{\underline{d}}$ and $V_{\underline{d}''}$ be $r+1$-dimensional subspaces of $H^{0}(\mathcal{L}_{\underline{d}})$ and $H^{0}(\mathcal{L}_{\underline{d}''})$, respectively, such that $\varphi_{\underline{d},\underline{d}''}(V_{\underline{d}})=V_{\underline{d}''}^{X_{1},0}$. Then 
\begin{center}
{\em dim}\,$\left(\dfrac{V_{\underline{d}}}{V_{\underline{d}}^{X_{2},0}+V_{\underline{d}}^{X_{3},0}}\right)=$\,{\em dim}\,$\left(\dfrac{V_{\underline{d}''}^{X_{1},0}+V_{\underline{d}''}^{X_{2},0}+V_{\underline{d}''}^{X_{3},0}}{V_{\underline{d}''}^{X_{2},0}+V_{\underline{d}''}^{X_{3},0}}\right)\,+$\,{\em dim}\,$\left(\dfrac{V_{\underline{d}''}^{X_{1},0}\cap(V_{\underline{d}''}^{X_{2},0}+V_{\underline{d}''}^{X_{3},0})}{V_{\underline{d}''}^{X_{3}^{c},0}\oplus V_{\underline{d}''}^{X_{2}^{c},0}}\right)$.
\end{center}
2. Let $V_{\underline{d}}$ and $V_{\underline{d}'''}$ be $r+1$-dimensional subspaces of $H^{0}(\mathcal{L}_{\underline{d}})$ and $H^{0}(\mathcal{L}_{\underline{d}'''})$, respectively, such that $\varphi_{\underline{d}''',\underline{d}}(V_{\underline{d}'''})=V_{\underline{d}}^{X_{3},0}$. Then
\begin{center}
{\em dim}\,$\left(\dfrac{V_{\underline{d}'''}}{V_{\underline{d}'''}^{X_{1},0}+V_{\underline{d}'''}^{X_{2},0}}\right)=$\,{\em dim}\,$\left(\dfrac{V_{\underline{d}}^{X_{1},0}+V_{\underline{d}}^{X_{2},0}+V_{\underline{d}}^{X_{3},0}}{V_{\underline{d}}^{X_{1},0}+V_{\underline{d}}^{X_{2},0}}\right)\,+$\,{\em dim}\,$\left(\dfrac{V_{\underline{d}}^{X_{3},0}\cap(V_{\underline{d}}^{X_{1},0}+V_{\underline{d}}^{X_{2},0})}{V_{\underline{d}}^{X_{2}^{c},0}\oplus V_{\underline{d}}^{X_{1}^{c},0}}\right)$.
\end{center}
\end{prop}
{\em Proof.} We will only prove the statement 1, as the proof of the statement 2 is analogous. Set $\beta:=\text{dim}\,V_{\underline{d}''}^{X_{1},0}-\text{dim}\,(V_{\underline{d}''}^{X_{3}^{c},0}\oplus V_{\underline{d}''}^{X_{2}^{c},0})$. Then, we have
\begin{equation}\label{dimension de cocientes1}
\beta=\text{dim}\,\left(\dfrac{V_{\underline{d}''}^{X_{1},0}}{V_{\underline{d}''}^{X_{1},0}\cap(V_{\underline{d}''}^{X_{2},0}+V_{\underline{d}''}^{X_{3},0})}\right)+\text{dim}\,\left(\dfrac{V_{\underline{d}''}^{X_{1},0}\cap(V_{\underline{d}''}^{X_{2},0}+V_{\underline{d}''}^{X_{3},0})}{V_{\underline{d}''}^{X_{3}^{c},0}\oplus V_{\underline{d}''}^{X_{2}^{c},0}}\right).
\end{equation}
On the other hand, since 
\begin{equation*}
\dfrac{V_{\underline{d}''}^{X_{1},0}}{V_{\underline{d}''}^{X_{1},0}\cap(V_{\underline{d}''}^{X_{2},0}+V_{\underline{d}''}^{X_{3},0})}\cong 
\dfrac{V_{\underline{d}''}^{X_{1},0}+V_{\underline{d}''}^{X_{2},0}+V_{\underline{d}''}^{X_{3},0}}{V_{\underline{d}''}^{X_{2},0}+V_{\underline{d}''}^{X_{3},0}},
\end{equation*}
we have
\begin{equation}\label{dimension de cocientes2}
\text{dim}\,\left(\dfrac{V_{\underline{d}''}^{X_{1},0}}{V_{\underline{d}''}^{X_{1},0}\cap(V_{\underline{d}''}^{X_{2},0}+V_{\underline{d}''}^{X_{3},0})}\right)=\text{dim}\,\left(\dfrac{V_{\underline{d}''}^{X_{1},0}+V_{\underline{d}''}^{X_{2},0}+V_{\underline{d}''}^{X_{3},0}}{V_{\underline{d}''}^{X_{2},0}+V_{\underline{d}''}^{X_{3},0}}\right).
\end{equation}
Also, by Proposition \ref{codimension1}, item 2, we have
\begin{equation}\label{dimension de cocientes3}
\beta=\text{dim}\,\left(\dfrac{V_{\underline{d}}}{V_{\underline{d}}^{X_{2},0}+V_{\underline{d}}^{X_{3},0}}\right).
\end{equation}
From (\ref{dimension de cocientes1}), (\ref{dimension de cocientes2}) and (\ref{dimension de cocientes3}), we get the desired equality.
\hfill $\Box$

We now introduce the notion of distributivity.
\begin{defn}\label{distributivity definition}
{\em Let $V_1, V_2$ and $V_3$ be vector subspaces of a $N$-dimensional vector space $V$. We will say that $V_1$ {\em distributes over} $V_2$ and $V_3$ if $V_1\cap (V_2+V_3)=V_1\cap V_2+V_1\cap V_3$.
}
\end{defn}

\begin{prop}\label{distributivity symmetry}
The notion of distributivity in Definition \ref{distributivity definition} is symmetric on the subspaces, that is, the following statements are equivalent:\\
1. $V_1$ distributes over $V_2$ and $V_3$.\\
2. $V_2$ distributes over $V_1$ and $V_3$.\\
3. $V_3$ distributes over $V_1$ and $V_2$.
\end{prop}
{\em Proof.} Muñoz [6], Remark 3.7.
\hfill $\Box$

\begin{defn}\label{distributivity in a subspace}
{\em Given $\underline{d}\geq \underline{0}$ and a subspace $V\subseteq H^{0}(\mathcal{L}_{\underline{d}})$, we will say that {\em the distributivity} holds in $V$ if $V^{X_{1},0}$ distributes over $V^{X_{2},0}$ and $V^{X_{3},0}$, that is,
\begin{center}
$V^{X_{1},0}\cap (V^{X_{2},0}+V^{X_{3},0})=V^{X_{1},0}\cap V^{X_{2},0}+V^{X_{1},0}\cap V^{X_{3},0}$.
\end{center}
}
\end{defn}

The following proposition states characterizations of distributivity under specific exactness conditions.

\begin{prop}\label{distributivity characterization}
The following statements hold:\\
1. Let $V_{\underline{d}}$ and $V_{\underline{d}''}$ be $r+1$-dimensional subspaces of $H^{0}(\mathcal{L}_{\underline{d}})$ and $H^{0}(\mathcal{L}_{\underline{d}''})$, respectively, such that $\varphi_{\underline{d},\underline{d}''}(V_{\underline{d}})=V_{\underline{d}''}^{X_{1},0}$. Then the distributivity holds in $V_{\underline{d}''}$ if and only if
\begin{center}
{\em dim}\,$(V_{\underline{d}}^{X_{2},0}+V_{\underline{d}}^{X_{3},0})\,-\,${\em dim}\,$(V_{\underline{d}''}^{X_{2},0}+V_{\underline{d}''}^{X_{3},0})=r+1\,-\,${\em dim}\,$(V_{\underline{d}''}^{X_{1},0}+V_{\underline{d}''}^{X_{2},0}+V_{\underline{d}''}^{X_{3},0})$
\end{center}
2. Let $V_{\underline{d}}$ and $V_{\underline{d}'''}$ be $r+1$-dimensional subspaces of $H^{0}(\mathcal{L}_{\underline{d}})$ and $H^{0}(\mathcal{L}_{\underline{d}'''})$, respectively, such that $\varphi_{\underline{d}''',\underline{d}}(V_{\underline{d}'''})=V_{\underline{d}}^{X_{3},0}$. Then the distributivity holds in $V_{\underline{d}}$ if and only if
\begin{center}
{\em dim}\,$(V_{\underline{d}'''}^{X_{1},0}+V_{\underline{d}'''}^{X_{2},0})\,-\,${\em dim}\,$(V_{\underline{d}}^{X_{1},0}+V_{\underline{d}}^{X_{2},0})=r+1\,-\,${\em dim}\,$(V_{\underline{d}}^{X_{1},0}+V_{\underline{d}}^{X_{2},0}+V_{\underline{d}}^{X_{3},0})$,
\end{center}
\end{prop}
{\em Proof.} We will only prove the statement 1, as the proof of the statement 2 is analogous because of the symmetry of the distributivity. Notice that the distributivity holds in $V_{\underline{d}''}$ if and only if 
\begin{center}
$V_{\underline{d}''}^{X_{1},0}\cap(V_{\underline{d}''}^{X_{2},0}+V_{\underline{d}''}^{X_{3},0})=V_{\underline{d}''}^{X_{3}^{c},0}\oplus V_{\underline{d}''}^{X_{2}^{c},0}$,
\end{center}
which is equivalent to 
\begin{equation*}
\text{dim}\left(\dfrac{V_{\underline{d}''}^{X_{1},0}\cap(V_{\underline{d}''}^{X_{2},0}+V_{\underline{d}''}^{X_{3},0})}{V_{\underline{d}''}^{X_{3}^{c},0}\oplus V_{\underline{d}''}^{X_{2}^{c},0}}\right)=0
\end{equation*}
On the other hand, by Proposition \ref{dimension de cocientes}, item 1, we have
\begin{equation*}
\text{dim}\left(\dfrac{V_{\underline{d}}}{V_{\underline{d}}^{X_{2},0}+V_{\underline{d}}^{X_{3},0}}\right)=\text{dim}\left(\dfrac{V_{\underline{d}''}^{X_{1},0}+V_{\underline{d}''}^{X_{2},0}+V_{\underline{d}''}^{X_{3},0}}{V_{\underline{d}''}^{X_{2},0}+V_{\underline{d}''}^{X_{3},0}}\right)+\,\text{dim}\left(\dfrac{V_{\underline{d}''}^{X_{1},0}\cap(V_{\underline{d}''}^{X_{2},0}+V_{\underline{d}''}^{X_{3},0})}{V_{\underline{d}''}^{X_{3}^{c},0}\oplus V_{\underline{d}''}^{X_{2}^{c},0}}\right)
\end{equation*}
Thus, the distributivity holds in $V_{\underline{d}''}$ if and only if 
\begin{equation*}
\text{dim}\left(\dfrac{V_{\underline{d}}}{V_{\underline{d}}^{X_{2},0}+V_{\underline{d}}^{X_{3},0}}\right)=\text{dim}\left(\dfrac{V_{\underline{d}''}^{X_{1},0}+V_{\underline{d}''}^{X_{2},0}+V_{\underline{d}''}^{X_{3},0}}{V_{\underline{d}''}^{X_{2},0}+V_{\underline{d}''}^{X_{3},0}}\right),
\end{equation*}
that is, 
\begin{center}
$r+1-\text{dim}\,(V_{\underline{d}}^{X_{2},0}+V_{\underline{d}}^{X_{3},0})=\text{dim}\,(V_{\underline{d}''}^{X_{1},0}+V_{\underline{d}''}^{X_{2},0}+V_{\underline{d}''}^{X_{3},0})-\text{dim}\,(V_{\underline{d}''}^{X_{2},0}+V_{\underline{d}''}^{X_{3},0})$.
\end{center}
From this, the desired result follows immediately.
\hfill $\Box$

\begin{prop}\label{suficiente1}
Let $V$ be a subspace of $H^{0}(\mathcal{L}_{\underline{d}})$ such that both $V^{X_{1},0}$ and $V^{X_{3},0}$ are contained in $V^{X_{2},0}$. Then the distributivity holds in $V$. 
\end{prop}
{\em Proof.} Since $V^{X_{1},0}\subseteq V^{X_{2},0}$ and $V^{X_{3},0}\subseteq V^{X_{2},0}$, we have $V^{X_{1},0}+V^{X_{3},0}\subseteq V^{X_{2},0}$, and hence
\begin{equation}\label{suficiente1ec1}
V^{X_{2},0}\cap(V^{X_{1},0}+V^{X_{3},0})=V^{X_{1},0}+V^{X_{3},0}.
\end{equation}
On the other hand, since $V^{X_{1},0}\subseteq V^{X_{2},0}$ and $V^{X_{3},0}\subseteq V^{X_{2},0}$, we also conclude that
\begin{equation}\label{suficiente1ec2}
V^{X_{1},0}=V^{X_{3}^{c},0}\,\, \text{and}\,\, V^{X_{3},0}=V^{X_{1}^{c},0}.
\end{equation}
It follows from (\ref{suficiente1ec1}) and (\ref{suficiente1ec2}) that
\begin{center}
$V^{X_{2},0}\cap(V^{X_{1},0}+V^{X_{3},0})=V^{X_{3}^{c},0}+V^{X_{1}^{c},0}$,
\end{center}
which proves the proposition.
\hfill $\Box$

\begin{prop}\label{series exactas desigualdad}
Let $\{(\mathcal{L}_{\underline{d}},V_{\underline{d}})\}_{\underline{d}}$ be an exact limit linear series of degree $d$ and dimension $r$. Then 
{\em
\begin{center}
$\displaystyle\sum_{\underline{d}}\text{dim}\,\left(\dfrac{V_{\underline{d}}}{V_{\underline{d}}^{X_{1},0}+V_{\underline{d}}^{X_{2},0}+V_{\underline{d}}^{X_{3},0}}\right)\geq r+1$
\end{center}
}
\end{prop}
{\em Proof.} It follows from Proposition \ref{dimension de cocientes}, item 1, that, for any $\underline{d}''$
\begin{equation*}
\text{dim}\left(\dfrac{V_{\underline{d}}}{V_{\underline{d}}^{X_{2},0}+V_{\underline{d}}^{X_{3},0}}\right)\geq \text{dim}\left(\dfrac{V_{\underline{d}''}^{X_{1},0}+V_{\underline{d}''}^{X_{2},0}+V_{\underline{d}''}^{X_{3},0}}{V_{\underline{d}''}^{X_{2},0}+V_{\underline{d}''}^{X_{3},0}}\right),
\end{equation*}
that is,
\begin{center}
$r+1-\text{dim}\,(V_{\underline{d}}^{X_{2},0}+V_{\underline{d}}^{X_{3},0})\geq\text{dim}\,(V_{\underline{d}''}^{X_{1},0}+V_{\underline{d}''}^{X_{2},0}+V_{\underline{d}''}^{X_{3},0})-\text{dim}\,(V_{\underline{d}''}^{X_{2},0}+V_{\underline{d}''}^{X_{3},0})$.
\end{center}
Then, for any $\underline{d}''$
{\em
\begin{center}
$r+1\,-\,${\em dim}\,$(V_{\underline{d}''}^{X_{1},0}+V_{\underline{d}''}^{X_{2},0}+V_{\underline{d}''}^{X_{3},0})\geq\,\,${\em dim}\,$(V_{\underline{d}}^{X_{2},0}+V_{\underline{d}}^{X_{3},0})\,-\,${\em dim}\,$(V_{\underline{d}''}^{X_{2},0}+V_{\underline{d}''}^{X_{3},0})$.
\end{center}
}
Thus, for each $0\leq l\leq d-1$, we have
\begin{align}
\displaystyle\sum_{i=1}^{d-l}\,(r+1-\text{dim}\,(V_{i-1\,l}^{X_{1},0}+V_{i-1\,l}^{X_{2},0}+V_{i-1\,l}^{X_{3},0}))\geq & \sum_{i=1}^{d-l}\,(\text{dim}\,(V_{il}^{X_{2},0}+V_{il}^{X_{3},0})-\text{dim}\,(V_{i-1\,l}^{X_{2},0}+V_{i-1\,l}^{X_{3},0})) \nonumber \\
=\,& \text{dim}\,(V_{d-l\,l}^{X_{2},0}+V_{d-l\,l}^{X_{3},0})-\text{dim}\,(V_{0l}^{X_{2},0}+V_{0l}^{X_{3},0}). \nonumber
\end{align}
It follows that
\begin{equation*}
\displaystyle\sum_{l=0}^{d-1}\sum_{i=1}^{d-l}\,(r+1-\text{dim}\,(V_{i-1\,l}^{X_{1},0}+V_{i-1\,l}^{X_{2},0}+V_{i-1\,l}^{X_{3},0}))\geq \sum_{l=0}^{d-1}\,\text{dim}\,(V_{d-l\,l}^{X_{2},0}+V_{d-l\,l}^{X_{3},0})-\sum_{l=0}^{d-1}\,\text{dim}\,(V_{0l}^{X_{2},0}+V_{0l}^{X_{3},0}),
\end{equation*}
that is,
\begin{equation*}
\displaystyle\sum_{i+l\leq d-1}\,(r+1-\text{dim}\,(V_{il}^{X_{1},0}+V_{il}^{X_{2},0}+V_{il}^{X_{3},0}))\geq \sum_{l=0}^{d-1}\,\text{dim}\,(V_{d-l\,l}^{X_{2},0}+V_{d-l\,l}^{X_{3},0})-\sum_{l=0}^{d-1}\,\text{dim}\,(V_{0l}^{X_{2},0}+V_{0l}^{X_{3},0}).
\end{equation*}
Then
\begin{align}\label{series exactas desigualdad1}
\displaystyle\sum_{i+l\leq d}\,(r+1-\text{dim}\,(V_{il}^{X_{1},0}+V_{il}^{X_{2},0}+V_{il}^{X_{3},0}))\geq &\sum_{l=0}^{d}\,(r+1-\text{dim}\,(V_{d-l\,l}^{X_{1},0}+V_{d-l\,l}^{X_{2},0}+V_{d-l\,l}^{X_{3},0})) \nonumber \\
&+\sum_{l=0}^{d-1}\,\text{dim}\,(V_{d-l\,l}^{X_{2},0}+V_{d-l\,l}^{X_{3},0}) \nonumber \\
&-\sum_{l=0}^{d-1}\,\text{dim}\,(V_{0l}^{X_{2},0}+V_{0l}^{X_{3},0}).
\end{align}
On the other hand, by Remark \ref{lado diagonal}, for each $0\leq l\leq d$, we have $V_{d-l\,l}^{X_{1},0}\subseteq V_{d-l\,l}^{X_{2},0}$, and hence
\begin{equation}\label{series exactas desigualdad2}
V_{d-l\,l}^{X_{2},0}+V_{d-l\,l}^{X_{3},0}=V_{d-l\,l}^{X_{1},0}+V_{d-l\,l}^{X_{2},0}+V_{d-l\,l}^{X_{3},0}.
\end{equation}
It follows from (\ref{series exactas desigualdad1}) and (\ref{series exactas desigualdad2}) that
\begin{align}
\displaystyle\sum_{i+l\leq d}\,(r+1-\text{dim}\,(V_{il}^{X_{1},0}+V_{il}^{X_{2},0}+V_{il}^{X_{3},0}))\geq &\sum_{l=0}^{d}\,(r+1-\text{dim}\,(V_{d-l\,l}^{X_{2},0}+V_{d-l\,l}^{X_{3},0})) \nonumber \\
&+\sum_{l=0}^{d-1}\,\text{dim}\,(V_{d-l\,l}^{X_{2},0}+V_{d-l\,l}^{X_{3},0}) \nonumber \\
&-\sum_{l=0}^{d-1}\,\text{dim}\,(V_{0l}^{X_{2},0}+V_{0l}^{X_{3},0}), \nonumber
\end{align}
that is,
\begin{equation}\label{series exactas desigualdad3}
\displaystyle\sum_{i+l\leq d}\,(r+1-\text{dim}\,(V_{il}^{X_{1},0}+V_{il}^{X_{2},0}+V_{il}^{X_{3},0}))\geq \sum_{l=0}^{d}\,(r+1-\text{dim}\,(V_{0l}^{X_{2},0}+V_{0l}^{X_{3},0})).
\end{equation}
On the other hand, by Corollary \ref{suma1corolario1}, for each $0\leq l\leq d-1$, we have
\begin{center}
$r+1-\text{dim}\,(V_{0l}^{X_{2},0}+V_{0l}^{X_{3},0})=\text{dim}\,V_{0\,l+1}^{X_{2},0}-\text{dim}\,V_{0l}^{X_{2},0}$.
\end{center}
Then
\begin{align}\label{series exactas desigualdad4}
\displaystyle\sum_{l=0}^{d}\,(r+1-\text{dim}\,(V_{0l}^{X_{2},0}+V_{0l}^{X_{3},0}))=&\,r+1-\text{dim}\,(V_{0d}^{X_{2},0}+V_{0d}^{X_{3},0}) \nonumber \\
&+\sum_{l=0}^{d-1}\,(\text{dim}\,V_{0\,l+1}^{X_{2},0}-\text{dim}\,V_{0l}^{X_{2},0}) \nonumber \\
=&\,r+1-\text{dim}\,(V_{0d}^{X_{2},0}+V_{0d}^{X_{3},0}) \nonumber \\
&+\text{dim}\,V_{0d}^{X_{2},0}-\text{dim}\,V_{00}^{X_{2},0}.
\end{align}
Now, from Remark \ref{lado diagonal}, we have
\begin{equation}\label{series exactas desigualdad5}
V_{0d}^{X_{2},0}+V_{0d}^{X_{3},0}=V_{0d}^{X_{2},0},
\end{equation}
and by Remark \ref{vertice superior derecho}, we have
\begin{equation}\label{series exactas desigualdad6}
V_{00}^{X_{2},0}=0.
\end{equation}
It follows from (\ref{series exactas desigualdad4}), (\ref{series exactas desigualdad5}) and (\ref{series exactas desigualdad6}) that
\begin{equation}\label{series exactas desigualdad7}
\displaystyle\sum_{l=0}^{d}\,(r+1-\text{dim}\,(V_{0l}^{X_{2},0}+V_{0l}^{X_{3},0}))=r+1.
\end{equation}
Thus, from (\ref{series exactas desigualdad3}) and (\ref{series exactas desigualdad7}), we obtain
\begin{center}
$\displaystyle\sum_{i+l\leq d}\,(r+1-\text{dim}\,(V_{il}^{X_{1},0}+V_{il}^{X_{2},0}+V_{il}^{X_{3},0}))\geq r+1$,
\end{center}
that is,
\begin{center}
$\displaystyle\sum_{\underline{d}}\text{dim}\,\left(\dfrac{V_{\underline{d}}}{V_{\underline{d}}^{X_{1},0}+V_{\underline{d}}^{X_{2},0}+V_{\underline{d}}^{X_{3},0}}\right)\geq r+1$.
\end{center}
\hfill $\Box$

We have the following corollary of the proof of Proposition \ref{series exactas desigualdad}.
\begin{cor}\label{sum of codimensions}
Let $\{(\mathcal{L}_{\underline{d}},V_{\underline{d}})\}_{\underline{d}}$ be an exact limit linear series of degree $d$ and dimension $r$. Then
{\em
\begin{center}
$\displaystyle\sum_{\underline{d}}\text{dim}\,\left(\dfrac{V_{\underline{d}}}{V_{\underline{d}}^{X_{1},0}+V_{\underline{d}}^{X_{2},0}+V_{\underline{d}}^{X_{3},0}}\right)=r+1$
\end{center}
}
if and only if the distributivity holds in $V_{\underline{d}}$ for any $\underline{d}$.
\end{cor}
{\em Proof.} It follows from the proof of Proposition \ref{series exactas desigualdad} that equality holds in 
\begin{center}
$\displaystyle\sum_{\underline{d}}\text{dim}\,\left(\dfrac{V_{\underline{d}}}{V_{\underline{d}}^{X_{1},0}+V_{\underline{d}}^{X_{2},0}+V_{\underline{d}}^{X_{3},0}}\right)\geq r+1$
\end{center}
if and only if equality holds in 
{\em
\begin{center}
$r+1\,-\,${\em dim}\,$(V_{\underline{d}''}^{X_{1},0}+V_{\underline{d}''}^{X_{2},0}+V_{\underline{d}''}^{X_{3},0})\geq\,\,${\em dim}\,$(V_{\underline{d}}^{X_{2},0}+V_{\underline{d}}^{X_{3},0})\,-\,${\em dim}\,$(V_{\underline{d}''}^{X_{2},0}+V_{\underline{d}''}^{X_{3},0})$,
\end{center}
}
for any $\underline{d}''$. Thus, by Proposition \ref{distributivity characterization}, item 1, we have that
\begin{center}
$\displaystyle\sum_{\underline{d}}\text{dim}\,\left(\dfrac{V_{\underline{d}}}{V_{\underline{d}}^{X_{1},0}+V_{\underline{d}}^{X_{2},0}+V_{\underline{d}}^{X_{3},0}}\right)=r+1$
\end{center}
if and only if the distributivity holds in $V_{\underline{d}''}$ for any $\underline{d}''$.

So to complete the proof, we will show that the distributivity always holds in $V_{il}$, if $i+l=d$. Let $i\geq 0$ and $l\geq 0$ such that $i+l=d$. By Remark \ref{lado diagonal}, we have that $V_{il}^{X_{1},0}\subseteq V_{il}^{X_{2},0}$ and $V_{il}^{X_{3},0}\subseteq V_{il}^{X_{2},0}$. Then, by Proposition \ref{suficiente1}, the distributivity holds in $V_{il}$, so we are done.
\hfill $\Box$

In the following result, we state that any simple limit linear series satisfies the equality in Corollary \ref{sum of codimensions}. We also characterize the multidegrees $\underline{d}^{1},\ldots,\underline{d}^{m}$ and the numbers $r_{1},\ldots,r_{m}$ in Definition \ref{simple lls definition}.

\begin{prop}\label{multidegrees characterization}
Let $\{(\mathcal{L}_{\underline{d}},V_{\underline{d}})\}_{\underline{d}}$ be a simple limit linear series of degree $d$ and dimension $r$. Then the following statements hold:\\
1. {\em $\displaystyle\sum_{\underline{d}}\text{dim}\,\left(\dfrac{V_{\underline{d}}}{V_{\underline{d}}^{X_{1},0}+V_{\underline{d}}^{X_{2},0}+V_{\underline{d}}^{X_{3},0}}\right)=r+1$. }\\
2. Let $\underline{d}^{1},\ldots,\underline{d}^{m}$ and $s^{\alpha}_{1},\ldots,s^{\alpha}_{r_{\alpha}}\in V_{\underline{d}^{\alpha}}$, for each $\alpha$, be as in Definition \ref{simple lls definition}. Then
\begin{center}
$V_{\underline{d}}^{X_{1},0}+V_{\underline{d}}^{X_{2},0}+V_{\underline{d}}^{X_{3},0}\neq V_{\underline{d}}$ if and only if $\underline{d}=\underline{d}^{\alpha}$ for some $\alpha$.
\end{center}
Furthermore, for each $\alpha$, we have
\begin{center}
$V_{\underline{d}^{\alpha}}=(V_{\underline{d}^{\alpha}}^{X_{1},0}+V_{\underline{d}^{\alpha}}^{X_{2},0}+V_{\underline{d}^{\alpha}}^{X_{3},0})\oplus \left\langle s^{\alpha}_{1},\ldots,s^{\alpha}_{r_{\alpha}}\right\rangle$.
\end{center}
\end{prop}
{\em Proof.} Let $\underline{d}^{1},\ldots,\underline{d}^{m}$ and $s^{\alpha}_{1},\ldots,s^{\alpha}_{r_{\alpha}}\in V_{\underline{d}^{\alpha}}$, for each $\alpha$, be as in Definition \ref{simple lls definition}. For each $\underline{d}^{\alpha}$, put $W_{\underline{d}^{\alpha}}:=\left\langle s^{\alpha}_{1},\ldots,s^{\alpha}_{r_{\alpha}}\right\rangle$, and for each $\underline{d}\neq \underline{d}^{\alpha}$ for any $\alpha$, put $W_{\underline{d}}:=0$. Then, for any $\underline{d}$, we have
\begin{equation}\label{multidegrees characterization0}
V_{\underline{d}}=\left(\bigoplus_{\underline{d}^{\alpha}\neq \underline{d}}\,\varphi_{\underline{d}^{\alpha},\underline{d}}(W_{\underline{d}^{\alpha}})\right)\oplus W_{\underline{d}}
\end{equation}
Thus,
\begin{align}
\displaystyle\sum_{\underline{d}}\,\text{dim}\,V_{\underline{d}}=&\sum_{\underline{d}}\,\left(\text{dim}\,\left(\bigoplus_{\underline{d}^{\alpha}\neq \underline{d}}\,\varphi_{\underline{d}^{\alpha},\underline{d}}(W_{\underline{d}^{\alpha}})\right)+\text{dim}\,W_{\underline{d}}\right) \nonumber \\
\leq&\sum_{\underline{d}}\,\left(\text{dim}\,(V_{\underline{d}}^{X_{1},0}+V_{\underline{d}}^{X_{2},0}+V_{\underline{d}}^{X_{3},0})+\text{dim}\,W_{\underline{d}}\right), \nonumber
\end{align}
where the inequality follows from the fact that
\begin{equation}\label{multidegrees characterization1}
\bigoplus_{\underline{d}^{\alpha}\neq \underline{d}}\,\varphi_{\underline{d}^{\alpha},\underline{d}}(W_{\underline{d}^{\alpha}})\subseteq V_{\underline{d}}^{X_{1},0}+V_{\underline{d}}^{X_{2},0}+V_{\underline{d}}^{X_{3},0},
\end{equation}
for any $\underline{d}$. Then
\begin{center}
$\displaystyle\sum_{\underline{d}}\,\text{dim}\,V_{\underline{d}}\leq \sum_{\underline{d}}\,\text{dim}\,(V_{\underline{d}}^{X_{1},0}+V_{\underline{d}}^{X_{2},0}+V_{\underline{d}}^{X_{3},0})+\sum_{\underline{d}}\,\text{dim}\,W_{\underline{d}}$,
\end{center}
and hence
\begin{center}
$\displaystyle\sum_{\underline{d}}\text{dim}\,\left(\dfrac{V_{\underline{d}}}{V_{\underline{d}}^{X_{1},0}+V_{\underline{d}}^{X_{2},0}+V_{\underline{d}}^{X_{3},0}}\right)\leq \sum_{\underline{d}}\,\text{dim}\,W_{\underline{d}}$.
\end{center}
Since $\displaystyle\sum_{\underline{d}}\,\text{dim}\,W_{\underline{d}}=r+1$, we get
\begin{equation}\label{multidegrees characterization2}
\displaystyle\sum_{\underline{d}}\text{dim}\,\left(\dfrac{V_{\underline{d}}}{V_{\underline{d}}^{X_{1},0}+V_{\underline{d}}^{X_{2},0}+V_{\underline{d}}^{X_{3},0}}\right)\leq r+1
\end{equation}
On the other hand, by Proposition \ref{simple is exact}, we have that $\{(\mathcal{L}_{\underline{d}},V_{\underline{d}})\}_{\underline{d}}$ is exact. Then, by Proposition \ref{series exactas desigualdad}, we have
\begin{equation}\label{multidegrees characterization3}
\displaystyle\sum_{\underline{d}}\text{dim}\,\left(\dfrac{V_{\underline{d}}}{V_{\underline{d}}^{X_{1},0}+V_{\underline{d}}^{X_{2},0}+V_{\underline{d}}^{X_{3},0}}\right)\geq r+1
\end{equation}
From (\ref{multidegrees characterization2}) and (\ref{multidegrees characterization3}), we get
\begin{center}
$\displaystyle\sum_{\underline{d}}\text{dim}\,\left(\dfrac{V_{\underline{d}}}{V_{\underline{d}}^{X_{1},0}+V_{\underline{d}}^{X_{2},0}+V_{\underline{d}}^{X_{3},0}}\right)=r+1$,
\end{center}
which proves the statement 1.

Now, since equality holds in (\ref{multidegrees characterization2}), we conclude that equality holds in (\ref{multidegrees characterization1}) for any $\underline{d}$, that is
\begin{equation}\label{multidegrees characterization4}
\bigoplus_{\underline{d}^{\alpha}\neq \underline{d}}\,\varphi_{\underline{d}^{\alpha},\underline{d}}(W_{\underline{d}^{\alpha}})=V_{\underline{d}}^{X_{1},0}+V_{\underline{d}}^{X_{2},0}+V_{\underline{d}}^{X_{3},0},
\end{equation}
for any $\underline{d}$. From (\ref{multidegrees characterization0}) and (\ref{multidegrees characterization4}) we obtain that
\begin{equation}\label{multidegrees characterization5}
V_{\underline{d}}=(V_{\underline{d}}^{X_{1},0}+V_{\underline{d}}^{X_{2},0}+V_{\underline{d}}^{X_{3},0})\oplus W_{\underline{d}},
\end{equation}
for any $\underline{d}$. Therefore, $V_{\underline{d}}^{X_{1},0}+V_{\underline{d}}^{X_{2},0}+V_{\underline{d}}^{X_{3},0}\neq V_{\underline{d}}$ if and only if $W_{\underline{d}}\neq 0$, which is equivalent to $\underline{d}=\underline{d}^{\alpha}$ for some $\alpha$.

Also, it follows from (\ref{multidegrees characterization5}) for $\underline{d}=\underline{d}^{\alpha}$ that
\begin{center}
$V_{\underline{d}^{\alpha}}=(V_{\underline{d}^{\alpha}}^{X_{1},0}+V_{\underline{d}^{\alpha}}^{X_{2},0}+V_{\underline{d}^{\alpha}}^{X_{3},0})\oplus \left\langle s^{\alpha}_{1},\ldots,s^{\alpha}_{r_{\alpha}}\right\rangle$,
\end{center}
for any $\alpha$. This proves the statement 2.
\hfill $\Box$

\begin{cor}\label{simple implica distributividad}
Let $\{(\mathcal{L}_{\underline{d}},V_{\underline{d}})\}_{\underline{d}}$ be a simple limit linear series of degree $d$ and dimension $r$. Then the distributivity holds in $V_{\underline{d}}$ for any $\underline{d}$.
\end{cor}
{\em Proof.} By Proposition \ref{multidegrees characterization}, item 1, we have that
\begin{center}
$\displaystyle\sum_{\underline{d}}\text{dim}\,\left(\dfrac{V_{\underline{d}}}{V_{\underline{d}}^{X_{1},0}+V_{\underline{d}}^{X_{2},0}+V_{\underline{d}}^{X_{3},0}}\right)=r+1$.
\end{center}
On the other hand, by Proposition \ref{simple is exact}, we have that $\{(\mathcal{L}_{\underline{d}},V_{\underline{d}})\}_{\underline{d}}$ is exact. Then, by Corollary \ref{sum of codimensions}, we conclude that the distributivity holds in $V_{\underline{d}}$ for any $\underline{d}$.
\hfill $\Box$

\begin{prop}\label{like simple basis}
Let $\{(\mathcal{L}_{\underline{d}},V_{\underline{d}})\}_{\underline{d}}$ be an exact limit linear series. Suppose that the distributivity holds in $V_{\underline{d}}$ for any $\underline{d}$. Then the following statements hold:\\
1. There exists a basis $\beta^{1}_{\underline{d}}$ for a subspace $W^{1}_{\underline{d}}\subseteq V_{\underline{d}}$ satisfying $V_{\underline{d}}=W^{1}_{\underline{d}}\oplus V_{\underline{d}}^{X_{1},0}$ for each $\underline{d}=(i,d-i-l,l)$ such that\\
(i) $\varphi_{\underline{d}''',\underline{d}}(\beta^{1}_{\underline{d}'''})\subseteq \beta^{1}_{\underline{d}}$ for any $i\geq 0$ and $l\geq 0$ such that $i+l\leq d-1$ (we say that the bases $\beta^{1}_{\underline{d}}$ increase in the vertical direction).\\
(ii) $\varphi_{\underline{d}',\underline{d}}(\beta^{1}_{\underline{d}'})\subseteq \beta^{1}_{\underline{d}}$ for any $i\geq 1$ and $l\geq 1$ such that $i+l\leq d$ (we say that the bases $\beta^{1}_{\underline{d}}$ increase in the diagonal direction).\\
(iii) $\varphi_{\underline{d}'',\underline{d}}(\beta^{1}_{\underline{d}''})\subseteq \beta^{1}_{\underline{d}}$ for any $i\geq 1$ and $l\geq 0$ such that $i+l\leq d$ (we say that the bases $\beta^{1}_{\underline{d}}$ increase in the horizontal direction).\\
2. There exists a basis $\beta^{2}_{\underline{d}}$ for a subspace $W^{2}_{\underline{d}}\subseteq V_{\underline{d}}$ satisfying $V_{\underline{d}}=W^{2}_{\underline{d}}\oplus V_{\underline{d}}^{X_{2},0}$ for each $\underline{d}=(i,d-i-l,l)$ such that\\
(i) $\varphi_{\underline{d}''',\underline{d}}(\beta^{2}_{\underline{d}'''})\subseteq \beta^{2}_{\underline{d}}$ for any $i\geq 0$ and $l\geq 0$ such that $i+l\leq d-1$ (we say that the bases $\beta^{2}_{\underline{d}}$ increase in the vertical direction).\\
(ii) $\varphi_{\underline{d},\underline{d}''}(\beta^{2}_{\underline{d}})\subseteq \beta^{2}_{\underline{d}''}$ for any $i\geq 1$ and $l\geq 0$ such that $i+l\leq d$ (we say that the bases $\beta^{2}_{\underline{d}}$ increase in the horizontal direction).\\
(iii) $\varphi_{\underline{d},\underline{d}'}(\beta^{2}_{\underline{d}})\subseteq \beta^{2}_{\underline{d}'}$ for any $i\geq 1$ and $l\geq 1$ such that $i+l\leq d$ (we say that the bases $\beta^{2}_{\underline{d}}$ increase in the diagonal direction).\\
3. There exists a basis $\beta^{3}_{\underline{d}}$ for a subspace $W^{3}_{\underline{d}}\subseteq V_{\underline{d}}$ satisfying $V_{\underline{d}}=W^{3}_{\underline{d}}\oplus V_{\underline{d}}^{X_{3},0}$ for each $\underline{d}=(i,d-i-l,l)$ such that\\
(i) $\varphi_{\underline{d},\underline{d}''}(\beta^{3}_{\underline{d}})\subseteq \beta^{3}_{\underline{d}''}$ for any $i\geq 1$ and $l\geq 0$ such that $i+l\leq d$ (we say that the bases $\beta^{3}_{\underline{d}}$ increase in the horizontal direction).\\
(ii) $\varphi_{\underline{d}',\underline{d}}(\beta^{3}_{\underline{d}'})\subseteq \beta^{3}_{\underline{d}}$ for any $i\geq 1$ and $l\geq 1$ such that $i+l\leq d$ (we say that the bases $\beta^{3}_{\underline{d}}$ increase in the diagonal direction).\\
(iii) $\varphi_{\underline{d},\underline{d}'''}(\beta^{3}_{\underline{d}})\subseteq \beta^{3}_{\underline{d}'''}$ for any $i\geq 0$ and $l\geq 0$ such that $i+l\leq d-1$ (we say that the bases $\beta^{3}_{\underline{d}}$ increase in the vertical direction).
\end{prop}
{\em Proof.} The proof of this proposition is essentially the proof given to Muñoz [5], Proposition 4.4. We will recall the construction of the subspaces $W^{q}_{\underline{d}}$ and the bases $\beta^{q}_{\underline{d}}$ (see Muñoz [5]). We will only see the cases $q=1$ and $q=2$, as the case $q=3$ is analogous to the case $q=1$. We will first see the case $q=1$.

We will use the notation $W^{1}_{il}:=W^{1}_{\underline{d}}$ and $\beta^{1}_{il}:=\beta^{1}_{\underline{d}}$ for each $\underline{d}=(i,d-i-l,l)$. The idea is first to construct the bases $\beta^{1}_{il}$ for $i=0$, then for $i=1$, and so on, until $i=d-1$. For each $i\geq 0$, we first construct the bases $\beta^{1}_{il}$ for $l=d-i$, then for $l=d-i-1$, and so on, until $l=0$.

For $i=0$, we inductively do the construction as follows:

{\em Step 1.} For $l=d$, we choose any linearly independent subset $\beta^{1}_{il}\subseteq V_{il}$ such that $V_{il}=\left\langle \beta^{1}_{il}\right\rangle\oplus V_{il}^{X_{1},0}$, and we set $W^{1}_{il}:=\left\langle \beta^{1}_{il}\right\rangle$.

{\em Step 2.} For $l=d-i-1,\ldots,l=0$, by induction, we have $V_{\underline{d}'''}=W^{1}_{\underline{d}'''}\oplus V_{\underline{d}'''}^{X_{1},0}$. This implies that $\varphi_{\underline{d}''',\underline{d}}(\beta^{1}_{\underline{d}'''})$ is a basis of $\varphi_{\underline{d}''',\underline{d}}(W^{1}_{\underline{d}'''})$ and $\varphi_{\underline{d}''',\underline{d}}(W^{1}_{\underline{d}'''})\cap V_{\underline{d}}^{X_{1},0}=0$. Thus, we can choose any linearly independent subset $\beta^{1}_{\underline{d}}\subseteq V_{\underline{d}}$ such that $\varphi_{\underline{d}''',\underline{d}}(\beta^{1}_{\underline{d}'''})\subseteq \beta^{1}_{\underline{d}}$ and $V_{\underline{d}}=\left\langle \beta^{1}_{\underline{d}}\right\rangle\oplus V_{\underline{d}}^{X_{1},0}$, and we set $W^{1}_{\underline{d}}:=\left\langle \beta^{1}_{\underline{d}}\right\rangle$.

Now, suppose inductively that, for a fixed $1\leq i\leq d-1$, the subspaces $W^{1}_{\tilde{i}l}\subseteq V_{\tilde{i}l}$ and the bases $\beta^{1}_{\tilde{i}l}$ have been constructed for $0\leq \tilde{i}\leq i-1$ and for $0\leq l\leq d-\tilde{i}$. Then, we inductively do the construction of the subspaces $W^{1}_{il}\subseteq V_{il}$ and the bases $\beta^{1}_{il}$ as follows:

{\em Step 1.} For $l=d-i$, we do the construction as in Step 2 above, with $\underline{d}'$ instead of $\underline{d}'''$. 

{\em Step 2.} For $l=d-i-1,\ldots,l=1$, by induction, we have 
\begin{center}
$V_{\underline{d}'}=W^{1}_{\underline{d}'}\oplus V_{\underline{d}'}^{X_{1},0}$\,,\, $V_{\underline{d}''}=W^{1}_{\underline{d}''}\oplus V_{\underline{d}''}^{X_{1},0}$\,,\, $V_{\underline{d}'''}=W^{1}_{\underline{d}'''}\oplus V_{\underline{d}'''}^{X_{1},0}$,\\
and $\varphi_{\underline{d}'',\underline{d}'}(\beta^{1}_{\underline{d}''})\subseteq \beta^{1}_{\underline{d}'}$, $\varphi_{\underline{d}'',\underline{d}'''}(\beta^{1}_{\underline{d}''})\subseteq \beta^{1}_{\underline{d}'''}$.
\end{center}
Then $\varphi_{\underline{d}',\underline{d}}(\beta^{1}_{\underline{d}'})$, $\varphi_{\underline{d}'',\underline{d}}(\beta^{1}_{\underline{d}''})$ and $\varphi_{\underline{d}''',\underline{d}}(\beta^{1}_{\underline{d}'''})$ are bases of the subspaces $\varphi_{\underline{d}',\underline{d}}(W^{1}_{\underline{d}'})$, $\varphi_{\underline{d}'',\underline{d}}(W^{1}_{\underline{d}''})$ and $\varphi_{\underline{d}''',\underline{d}}(W^{1}_{\underline{d}'''})$, respectively. Also, $\varphi_{\underline{d}'',\underline{d}}(W^{1}_{\underline{d}''})=V_{\underline{d}}^{X_{1}^{c},0}$. Furthermore, in Muñoz [5] it is shown that 
\begin{center}
$\varphi_{\underline{d}',\underline{d}}(W^{1}_{\underline{d}'})\cap \varphi_{\underline{d}''',\underline{d}}(W^{1}_{\underline{d}'''})=V_{\underline{d}}^{X_{1}^{c},0}$,
\end{center}
\begin{equation}\label{like simple basis1}
\varphi_{\underline{d}',\underline{d}}(\beta^{1}_{\underline{d}'})\cup \varphi_{\underline{d}''',\underline{d}}(\beta^{1}_{\underline{d}'''})\, \text{is a basis of}\, \varphi_{\underline{d}',\underline{d}}(W^{1}_{\underline{d}'})+\varphi_{\underline{d}''',\underline{d}}(W^{1}_{\underline{d}'''})
\end{equation}
\begin{equation}\label{like simple basis2}
\text{and}\, V_{\underline{d}}^{X_{2},0}+V_{\underline{d}}^{X_{3},0}=(V_{\underline{d}}^{X_{3}^{c},0}\oplus V_{\underline{d}}^{X_{2}^{c},0})\oplus (\varphi_{\underline{d}',\underline{d}}(W^{1}_{\underline{d}'})+\varphi_{\underline{d}''',\underline{d}}(W^{1}_{\underline{d}'''})).
\end{equation}
Since the distributivity holds in $V_{\underline{d}}$, we have that
\begin{equation}\label{like simple basis3}
V_{\underline{d}}^{X_{1},0}\cap (V_{\underline{d}}^{X_{2},0}+V_{\underline{d}}^{X_{3},0})=V_{\underline{d}}^{X_{3}^{c},0}\oplus V_{\underline{d}}^{X_{2}^{c},0}.
\end{equation}
It follows from (\ref{like simple basis2}) and (\ref{like simple basis3}) that 
\begin{equation}\label{like simple basis4}
(\varphi_{\underline{d}',\underline{d}}(W^{1}_{\underline{d}'})+\varphi_{\underline{d}''',\underline{d}}(W^{1}_{\underline{d}'''}))\cap V_{\underline{d}}^{X_{1},0}=0.
\end{equation}
From (\ref{like simple basis1}) and (\ref{like simple basis4}), we conclude that we can choose any linearly independent subset $\beta^{1}_{\underline{d}}\subseteq V_{\underline{d}}$ such that $\varphi_{\underline{d}',\underline{d}}(\beta^{1}_{\underline{d}'})\cup \varphi_{\underline{d}''',\underline{d}}(\beta^{1}_{\underline{d}'''})\subseteq \beta^{1}_{\underline{d}}$ and $V_{\underline{d}}=\left\langle \beta^{1}_{\underline{d}}\right\rangle\oplus V_{\underline{d}}^{X_{1},0}$, and we set $W^{1}_{\underline{d}}:=\left\langle \beta^{1}_{\underline{d}}\right\rangle$.

{\em Step 3.} For $l=0$, we do the construction as in Step 2 of the case $i=0$.

Finally, for $i=d$ and $l=0$, we do the construction as in Step 2 of the case $i=0$, with $\underline{d}''$ instead of $\underline{d}'''$. This finishes the construction for the case $q=1$.

Now, we will see the case $q=2$.

We will use the notation $W^{2}_{il}:=W^{2}_{\underline{d}}$ and $\beta^{2}_{il}:=\beta^{2}_{\underline{d}}$ for each $\underline{d}=(i,d-i-l,l)$. The idea is first to construct the bases $\beta^{2}_{il}$ for $i+l=d$, then for $i+l=d-1$, and so on, until $i+l=0$.

Suppose inductively that, for a fixed $2\leq m\leq d$, the subspaces $W^{2}_{il}\subseteq V_{il}$ and the bases $\beta^{2}_{il}$ have been constructed for $i+l\geq d-(m-1)$. Then, we inductively do the construction of the subspaces $W^{2}_{il}\subseteq V_{il}$ and the bases $\beta^{2}_{il}$, for each multidegree $(i,d-i-l,l)$ such that $i+l=d-m$, as follows:

Repeat Step 2 of the case $1\leq i\leq d-1$ of the constuction for the case $q=1$, with $(i+1,d-(i+1)-l,l)$, $(i+1,d-(i+1)-(l+1),l+1)$ and $(i,d-i-(l+1),l+1)$ instead of $\underline{d}'$, $\underline{d}''$ and $\underline{d}'''$, respectively, putting the superscript $q=2$ instead of $q=1$, and exchanging $X_{1}$ with $X_{2}$.

To complete the proof, we do the construction of the subspaces $W^{2}_{il}\subseteq V_{il}$ and the bases $\beta^{2}_{il}$, for each multidegree $(i,d-i-l,l)$ such that $i+l\geq d-1$, as follows:

For each $(i,d-i-l,l)$ such that $i+l=d$, we choose any linearly independent subset $\beta^{2}_{il}\subseteq V_{il}$ such that $V_{il}=\left\langle \beta^{2}_{il}\right\rangle\oplus V_{il}^{X_{2},0}$, and we set $W^{2}_{il}:=\left\langle \beta^{2}_{il}\right\rangle$.

For each multidegree $\underline{d}=(i,d-i-l,l)$ such that $i+l=d-1$, we construct the subspace $W^{2}_{il}\subseteq V_{il}$ and the basis $\beta^{2}_{il}$ by using an analogous reasoning to that of the above construction for $i+l=d-m$ with $2\leq m\leq d$, but without considering the multidegree $(i+1,d-(i+1)-(l+1),l+1)$, as in the present case, we have
\begin{center}
$d-(i+1)-(l+1)=d-(i+l)-2=d-(d-1)-2=-1<0$.
\end{center}
This finishes the construction for the case $q=2$.
\hfill $\Box$

\begin{prop}\label{prop1structure}
Let $\{(\mathcal{L}_{\underline{d}},V_{\underline{d}})\}_{\underline{d}}$ be an exact limit linear series, and suppose that the distributivity holds in $V_{\underline{d}}$ for any $\underline{d}$. For each $\underline{d}=(i,d-i-l,l)$, let $W^{1}_{\underline{d}}, W^{2}_{\underline{d}}$ and $W^{3}_{\underline{d}}$ be the subspaces of $V_{\underline{d}}$ constructed in the proof of Proposition \ref{like simple basis}, and set $\widetilde{\underline{d}}:=(i+1,d-i-l-1,l)$ for $i+l\leq d-1$. Then the following statements hold:\\
1. For $i=0$ and $l=d$, we have
\begin{center}
$V_{\underline{d}}^{X_{1},0}+V_{\underline{d}}^{X_{2},0}+V_{\underline{d}}^{X_{3},0}=\varphi_{\widetilde{\underline{d}}_{3},\underline{d}}(W^{3}_{\widetilde{\underline{d}}_{3}})$,
\end{center}
where $\widetilde{\underline{d}}_{3}=(i,d-i-l+1,l-1)$.\\
2. For $i=0$ and $0\leq l\leq d-1$, we have
\begin{center}
$V_{\underline{d}}^{X_{1},0}+V_{\underline{d}}^{X_{2},0}+V_{\underline{d}}^{X_{3},0}=\varphi_{\underline{d}''',\underline{d}}(W^{1}_{\underline{d}'''})+\varphi_{\underline{d}''',\underline{d}}(W^{2}_{\underline{d}'''})+\varphi_{\widetilde{\underline{d}},\underline{d}}(W^{2}_{\widetilde{\underline{d}}})+\varphi_{\widetilde{\underline{d}},\underline{d}}(W^{3}_{\widetilde{\underline{d}}})$.
\end{center}
3. For $1\leq i\leq d-1$ and $l=d-i$, we have
\begin{center}
$V_{\underline{d}}^{X_{1},0}+V_{\underline{d}}^{X_{2},0}+V_{\underline{d}}^{X_{3},0}=\varphi_{\underline{d}',\underline{d}}(W^{1}_{\underline{d}'})+\varphi_{\underline{d}',\underline{d}}(W^{3}_{\underline{d}'})$.
\end{center}
4. For $1\leq i\leq d-1$ and $1\leq l\leq d-i-1$, we have
\begin{align}
V_{\underline{d}}^{X_{1},0}+V_{\underline{d}}^{X_{2},0}+V_{\underline{d}}^{X_{3},0}=&\,\varphi_{\underline{d}',\underline{d}}(W^{1}_{\underline{d}'})+\varphi_{\underline{d}',\underline{d}}(W^{3}_{\underline{d}'}) \nonumber \\
&+\varphi_{\underline{d}''',\underline{d}}(W^{1}_{\underline{d}'''})+\varphi_{\underline{d}''',\underline{d}}(W^{2}_{\underline{d}'''})+\varphi_{\widetilde{\underline{d}},\underline{d}}(W^{2}_{\widetilde{\underline{d}}})+\varphi_{\widetilde{\underline{d}},\underline{d}}(W^{3}_{\widetilde{\underline{d}}}). \nonumber
\end{align}
5. For $1\leq i\leq d-1$ and $l=0$, we have
\begin{center}
$V_{\underline{d}}^{X_{1},0}+V_{\underline{d}}^{X_{2},0}+V_{\underline{d}}^{X_{3},0}=\varphi_{\underline{d}''',\underline{d}}(W^{1}_{\underline{d}'''})+\varphi_{\underline{d}''',\underline{d}}(W^{2}_{\underline{d}'''})+\varphi_{\widetilde{\underline{d}},\underline{d}}(W^{2}_{\widetilde{\underline{d}}})+\varphi_{\widetilde{\underline{d}},\underline{d}}(W^{3}_{\widetilde{\underline{d}}})$.
\end{center}
6. For $i=d$ and $l=0$, we have
\begin{center}
$V_{\underline{d}}^{X_{1},0}+V_{\underline{d}}^{X_{2},0}+V_{\underline{d}}^{X_{3},0}=\varphi_{\underline{d}'',\underline{d}}(W^{1}_{\underline{d}''})$.
\end{center}
\end{prop}
{\em Proof.} We will first prove the statement 1. For $i=0$ and $l=d$, we have $\underline{d}=(0,0,d)$. Then $V_{\underline{d}}^{X_{3},0}=0$ and $V_{\underline{d}}^{X_{1},0}=V_{\underline{d}}^{X_{2},0}$. Thus
\begin{equation}\label{prop1structure1}
V_{\underline{d}}^{X_{1},0}+V_{\underline{d}}^{X_{2},0}+V_{\underline{d}}^{X_{3},0}=V_{\underline{d}}^{X_{1},0}+V_{\underline{d}}^{X_{2},0}=V_{\underline{d}}^{X_{3}^{c},0}.
\end{equation}
On the other hand, since
\begin{center}
$V_{\widetilde{\underline{d}}_{3}}=W^{3}_{\widetilde{\underline{d}}_{3}}\oplus V_{\widetilde{\underline{d}}_{3}}^{X_{3},0}$,
\end{center}
we have
\begin{equation}\label{prop1structure2}
V_{\underline{d}}^{X_{3}^{c},0}=\varphi_{\widetilde{\underline{d}}_{3},\underline{d}}(V_{\widetilde{\underline{d}}_{3}})=\varphi_{\widetilde{\underline{d}}_{3},\underline{d}}(W^{3}_{\widetilde{\underline{d}}_{3}}).
\end{equation}
It follows from (\ref{prop1structure1}) and (\ref{prop1structure2}) that
\begin{center}
$V_{\underline{d}}^{X_{1},0}+V_{\underline{d}}^{X_{2},0}+V_{\underline{d}}^{X_{3},0}=\varphi_{\widetilde{\underline{d}}_{3},\underline{d}}(W^{3}_{\widetilde{\underline{d}}_{3}})$,
\end{center}
which proves the statement 1. 

Now, we will prove the statement 2. For $i=0$ and $0\leq l\leq d-1$, we have $\underline{d}=(0,d-l,l)$. Then $V_{\underline{d}}^{X_{2},0}\subseteq V_{\underline{d}}^{X_{1},0}$. Thus
\begin{equation}\label{prop1structure3}
V_{\underline{d}}^{X_{1},0}+V_{\underline{d}}^{X_{2},0}+V_{\underline{d}}^{X_{3},0}=V_{\underline{d}}^{X_{1},0}+V_{\underline{d}}^{X_{3},0}.
\end{equation}
Now, since 
\begin{center}
$V_{\underline{d}'''}=W^{1}_{\underline{d}'''}\oplus V_{\underline{d}'''}^{X_{1},0}$ and $V_{\underline{d}'''}^{X_{1},0}\supseteq V_{\underline{d}'''}^{X_{3}^{c},0}=\text{Ker}\,(\varphi_{\underline{d}''',\underline{d}}|_{V_{\underline{d}'''}})$,
\end{center}
we have
\begin{equation}\label{prop1structure4}
V_{\underline{d}}^{X_{3},0}=\varphi_{\underline{d}''',\underline{d}}(V_{\underline{d}'''})=\varphi_{\underline{d}''',\underline{d}}(W^{1}_{\underline{d}'''})\oplus V_{\underline{d}}^{X_{2}^{c},0}.
\end{equation}
Analogously, we have
\begin{equation}\label{prop1structure5}
V_{\underline{d}}^{X_{1},0}=\varphi_{\widetilde{\underline{d}},\underline{d}}(W^{3}_{\widetilde{\underline{d}}})\oplus V_{\underline{d}}^{X_{2}^{c},0}.
\end{equation}
It follows from (\ref{prop1structure4}) and (\ref{prop1structure5}) that
\begin{equation}\label{prop1structure6}
V_{\underline{d}}^{X_{1},0}+V_{\underline{d}}^{X_{3},0}=\varphi_{\underline{d}''',\underline{d}}(W^{1}_{\underline{d}'''})\oplus \varphi_{\widetilde{\underline{d}},\underline{d}}(W^{3}_{\widetilde{\underline{d}}})\oplus V_{\underline{d}}^{X_{2}^{c},0}.
\end{equation}
On the other hand, it follows from the proof of Proposition \ref{like simple basis} that
\begin{equation}\label{prop1structure7}
V_{\underline{d}}^{X_{2}^{c},0}=\varphi_{\widetilde{\underline{d}},\underline{d}}(W^{2}_{\widetilde{\underline{d}}})\cap \varphi_{\underline{d}''',\underline{d}}(W^{2}_{\underline{d}'''}).
\end{equation}
From (\ref{prop1structure6}) and (\ref{prop1structure7}), we obtain that
\begin{center}
$V_{\underline{d}}^{X_{1},0}+V_{\underline{d}}^{X_{3},0}=\varphi_{\underline{d}''',\underline{d}}(W^{1}_{\underline{d}'''})\oplus \varphi_{\widetilde{\underline{d}},\underline{d}}(W^{3}_{\widetilde{\underline{d}}})\oplus (\varphi_{\widetilde{\underline{d}},\underline{d}}(W^{2}_{\widetilde{\underline{d}}})\cap \varphi_{\underline{d}''',\underline{d}}(W^{2}_{\underline{d}'''}))$,
\end{center}
which implies
\begin{equation}\label{prop1structure8}
V_{\underline{d}}^{X_{1},0}+V_{\underline{d}}^{X_{3},0}\subseteq \varphi_{\underline{d}''',\underline{d}}(W^{1}_{\underline{d}'''})+\varphi_{\underline{d}''',\underline{d}}(W^{2}_{\underline{d}'''})+\varphi_{\widetilde{\underline{d}},\underline{d}}(W^{2}_{\widetilde{\underline{d}}})+ \varphi_{\widetilde{\underline{d}},\underline{d}}(W^{3}_{\widetilde{\underline{d}}}).
\end{equation}
It follows from (\ref{prop1structure3}) and (\ref{prop1structure8}) that
\begin{equation}\label{prop1structureasterisco}
V_{\underline{d}}^{X_{1},0}+V_{\underline{d}}^{X_{2},0}+V_{\underline{d}}^{X_{3},0}\subseteq \varphi_{\underline{d}''',\underline{d}}(W^{1}_{\underline{d}'''})+\varphi_{\underline{d}''',\underline{d}}(W^{2}_{\underline{d}'''})+\varphi_{\widetilde{\underline{d}},\underline{d}}(W^{2}_{\widetilde{\underline{d}}})+ \varphi_{\widetilde{\underline{d}},\underline{d}}(W^{3}_{\widetilde{\underline{d}}}).
\end{equation}
Since the other inclusion is obvious, equality holds in (\ref{prop1structureasterisco}), which proves the statement 2.

Now, we will prove the statement 3. For $1\leq i\leq d-1$ and $l=d-i$, we have $\underline{d}=(i,0,d-i)$. Then $V_{\underline{d}}^{X_{1},0}\subseteq V_{\underline{d}}^{X_{2},0}$ and  $V_{\underline{d}}^{X_{3},0}\subseteq V_{\underline{d}}^{X_{2},0}$. Thus
\begin{equation}\label{prop1structure9}
V_{\underline{d}}^{X_{1},0}+V_{\underline{d}}^{X_{2},0}+V_{\underline{d}}^{X_{3},0}=V_{\underline{d}}^{X_{2},0}.
\end{equation}
Now, since 
\begin{center}
$V_{\underline{d}'}=W^{1}_{\underline{d}'}\oplus V_{\underline{d}'}^{X_{1},0}$ and $V_{\underline{d}'}^{X_{1},0}\supseteq V_{\underline{d}'}^{X_{2}^{c},0}=\text{Ker}\,(\varphi_{\underline{d}',\underline{d}}|_{V_{\underline{d}'}})$,
\end{center}
we have
\begin{equation}\label{prop1structure10}
V_{\underline{d}}^{X_{2},0}=\varphi_{\underline{d}',\underline{d}}(V_{\underline{d}'})=\varphi_{\underline{d}',\underline{d}}(W^{1}_{\underline{d}'})\oplus V_{\underline{d}}^{X_{3}^{c},0}.
\end{equation}
On the other hand, as above, we have
\begin{equation}\label{prop1structure11}
V_{\underline{d}}^{X_{3}^{c},0}=\varphi_{\widetilde{\underline{d}}_{3},\underline{d}}(W^{3}_{\widetilde{\underline{d}}_{3}}),
\end{equation}
where $\widetilde{\underline{d}}_{3}=(i,d-i-l+1,l-1)$, as the same proof of (\ref{prop1structure2}) works here. From (\ref{prop1structure10}) and (\ref{prop1structure11}), we conclude that
\begin{equation}\label{prop1structure12}
V_{\underline{d}}^{X_{2},0}=\varphi_{\underline{d}',\underline{d}}(W^{1}_{\underline{d}'})\oplus \varphi_{\widetilde{\underline{d}}_{3},\underline{d}}(W^{3}_{\widetilde{\underline{d}}_{3}}).
\end{equation}
Now, from Proposition \ref{like simple basis}, item 3, we get $\varphi_{\widetilde{\underline{d}}_{3},\underline{d}'}(W^{3}_{\widetilde{\underline{d}}_{3}})\subseteq W^{3}_{\underline{d}'}$, so 
\begin{equation}\label{prop1structure13}
\varphi_{\widetilde{\underline{d}}_{3},\underline{d}}(W^{3}_{\widetilde{\underline{d}}_{3}})=\varphi_{\underline{d}',\underline{d}}(\varphi_{\widetilde{\underline{d}}_{3},\underline{d}'}(W^{3}_{\widetilde{\underline{d}}_{3}}))\subseteq \varphi_{\underline{d}',\underline{d}}(W^{3}_{\underline{d}'}).
\end{equation}
It follows from (\ref{prop1structure12}) and (\ref{prop1structure13}) that 
\begin{equation}\label{prop1structure14}
V_{\underline{d}}^{X_{2},0}\subseteq\varphi_{\underline{d}',\underline{d}}(W^{1}_{\underline{d}'})+\varphi_{\underline{d}',\underline{d}}(W^{3}_{\underline{d}'}).
\end{equation}
Then, from (\ref{prop1structure9}) and (\ref{prop1structure14}), we conclude that
\begin{equation}\label{prop1structure15}
V_{\underline{d}}^{X_{1},0}+V_{\underline{d}}^{X_{2},0}+V_{\underline{d}}^{X_{3},0}\subseteq\varphi_{\underline{d}',\underline{d}}(W^{1}_{\underline{d}'})+\varphi_{\underline{d}',\underline{d}}(W^{3}_{\underline{d}'}).
\end{equation}
Since the other inclusion is obvious, equality holds in (\ref{prop1structure15}), which proves the statement 3.

Now, we will prove the statement 4. Let $i,l$ be positive integers such that $i+l\leq d-1$. As above, we have
\begin{equation}\label{prop1structure16}
V_{\underline{d}}^{X_{1},0}+V_{\underline{d}}^{X_{3},0}\subseteq \varphi_{\underline{d}''',\underline{d}}(W^{1}_{\underline{d}'''})+\varphi_{\underline{d}''',\underline{d}}(W^{2}_{\underline{d}'''})+\varphi_{\widetilde{\underline{d}},\underline{d}}(W^{2}_{\widetilde{\underline{d}}})+ \varphi_{\widetilde{\underline{d}},\underline{d}}(W^{3}_{\widetilde{\underline{d}}})
\end{equation}
and
\begin{equation}\label{prop1structure17}
V_{\underline{d}}^{X_{2},0}\subseteq \varphi_{\underline{d}',\underline{d}}(W^{1}_{\underline{d}'})+\varphi_{\underline{d}',\underline{d}}(W^{3}_{\underline{d}'}),
\end{equation}
as the same proofs of (\ref{prop1structure8}) and (\ref{prop1structure14}) work here. It follows from (\ref{prop1structure16}) and (\ref{prop1structure17}) that
\begin{align}\label{prop1structure18}
V_{\underline{d}}^{X_{1},0}+V_{\underline{d}}^{X_{2},0}+V_{\underline{d}}^{X_{3},0}\subseteq & \,\varphi_{\underline{d}',\underline{d}}(W^{1}_{\underline{d}'})+\varphi_{\underline{d}',\underline{d}}(W^{3}_{\underline{d}'}) \nonumber \\
& +\varphi_{\underline{d}''',\underline{d}}(W^{1}_{\underline{d}'''})+\varphi_{\underline{d}''',\underline{d}}(W^{2}_{\underline{d}'''})+\varphi_{\widetilde{\underline{d}},\underline{d}}(W^{2}_{\widetilde{\underline{d}}})+ \varphi_{\widetilde{\underline{d}},\underline{d}}(W^{3}_{\widetilde{\underline{d}}}).
\end{align}
Since the other inclusion is obvious, equality holds in (\ref{prop1structure18}). This proves the statement 4.

Now, we will prove the statement 5. For $1\leq i\leq d-1$ and $l=0$, we have $\underline{d}=(i,d-i,0)$. Then $V_{\underline{d}}^{X_{2},0}\subseteq V_{\underline{d}}^{X_{3},0}$. Thus
\begin{center}
$V_{\underline{d}}^{X_{1},0}+V_{\underline{d}}^{X_{2},0}+V_{\underline{d}}^{X_{3},0}=V_{\underline{d}}^{X_{1},0}+V_{\underline{d}}^{X_{3},0}$.
\end{center}
From here, the proof follows as in the proof of the statement 2.

Finally, the proof of the statement 6 is analogous to that of the statement 1.
\hfill $\Box$

Before stating the next result, we introduce some notation. Let $M$ be the set of nonnegative multidegrees of total degree $d$. For each $\underline{d}=(i,d-i-l,l)$, we set
\begin{center}
$M^{1}_{\underline{d}}:=\{\underline{d}_{0}=(\tilde{i},d-\tilde{i}-\tilde{l},\tilde{l})\in M\,/\,\,\tilde{i}\leq i\,\,\text{and}\,\,\tilde{i}-i\leq \tilde{l}-l\}$, \\
$M^{2}_{\underline{d}}:=\{\underline{d}_{0}=(\tilde{i},d-\tilde{i}-\tilde{l},\tilde{l})\in M\,/\,\,\tilde{i}\geq i\,\,\text{and}\,\,\tilde{l}\geq l\}$ and \\
$M^{3}_{\underline{d}}:=\{\underline{d}_{0}=(\tilde{i},d-\tilde{i}-\tilde{l},\tilde{l})\in M\,/\,\,\tilde{l}\leq l\,\,\text{and}\,\,\tilde{l}-l\leq \tilde{i}-i\}$.
\end{center}
Note that $M^{1}_{\underline{d}}\cup M^{2}_{\underline{d}}\cup M^{3}_{\underline{d}}=M$.

\begin{prop}\label{prop2structure}
Let $\{(\mathcal{L}_{\underline{d}},V_{\underline{d}})\}_{\underline{d}}$ be an exact limit linear series, and suppose that the distributivity holds in $V_{\underline{d}}$ for any $\underline{d}$. Write
\begin{center}
$\{\underline{d}\,/\,\,V_{\underline{d}}^{X_{1},0}+V_{\underline{d}}^{X_{2},0}+V_{\underline{d}}^{X_{3},0}\neq V_{\underline{d}}\}=\{\underline{d}^{1},\ldots,\underline{d}^{m}\}$,
\end{center}
and for each $\alpha=1,\ldots,m$, let $\{s^{\alpha}_{1},\ldots,s^{\alpha}_{r_{\alpha}}\}\subseteq V_{\underline{d}^{\alpha}}$ be any linearly independent subset such that
\begin{center}
$V_{\underline{d}^{\alpha}}=(V_{\underline{d}}^{X_{1},0}+V_{\underline{d}}^{X_{2},0}+V_{\underline{d}}^{X_{3},0})\oplus \left\langle s^{\alpha}_{1},\ldots,s^{\alpha}_{r_{\alpha}}\right\rangle$.
\end{center}
Then, for each $q=1,2,3$ and for each $\underline{d}=(i,d-i-l,l)$, the basis $\beta^{q}_{\underline{d}}$ constructed in the proof of Proposition \ref{like simple basis} can be chosen as
\begin{center}
$\beta^{q}_{\underline{d}}=\{\varphi_{\underline{d}^{\alpha},\underline{d}}(s^{\alpha}_{z})\,/\,\,\underline{d}^{\alpha}\in M^{q}_{\underline{d}}\,\,\text{and}\,\,z=1,\ldots,r_{\alpha}\}$.
\end{center}
\end{prop}
{\em Proof.} We will only prove the cases $q=1$ and $q=2$, as the case $q=3$ is analogous to the case $q=1$. We will first prove the case $q=1$. 

For $i=0$, we prove the stated result inductively as follows:

For $l=d$, we have $\underline{d}=(0,0,d)$, so $V_{\underline{d}}^{X_{3},0}=0$ and $V_{\underline{d}}^{X_{1},0}=V_{\underline{d}}^{X_{2},0}$. In particular, 
\begin{center}
$V_{\underline{d}}^{X_{1},0}+V_{\underline{d}}^{X_{2},0}+V_{\underline{d}}^{X_{3},0}=V_{\underline{d}}^{X_{1},0}$. 
\end{center}
Then, from the proof of Proposition \ref{like simple basis}, we conclude that $\beta^{1}_{\underline{d}}$ can be chosen as
\begin{equation}\label{prop2structureasterisco}
\beta^{1}_{\underline{d}}=\{s^{\widetilde{\alpha}}_{1},\ldots,s^{\widetilde{\alpha}}_{r_{\widetilde{\alpha}}}\}\,\,\text{if}\,\,\underline{d}=\underline{d}^{\widetilde{\alpha}}\,\,\text{for some}\,\,\widetilde{\alpha},
\end{equation}
\begin{equation}\label{prop2structure1}
\text{and}\,\,\beta^{1}_{\underline{d}}=\emptyset\,\,\text{otherwise}.
\end{equation}
On the other hand, notice that 
\begin{equation}\label{prop2structureasteriscoprima}
M^{1}_{\underline{d}}=\{\underline{d}\}, 
\end{equation}
as $\underline{d}=(0,0,d)$. Thus, if $\underline{d}=\underline{d}^{\widetilde{\alpha}}$ for some $\widetilde{\alpha}$, we have
\begin{equation}\label{prop2structure2}
\{s^{\widetilde{\alpha}}_{1},\ldots,s^{\widetilde{\alpha}}_{r_{\widetilde{\alpha}}}\}=\{\varphi_{\underline{d}^{\alpha},\underline{d}}(s^{\alpha}_{z})\,/\,\,\underline{d}^{\alpha}\in M^{1}_{\underline{d}}\,\,\text{and}\,\,z=1,\ldots,r_{\alpha}\}.
\end{equation}
From (\ref{prop2structureasterisco}), (\ref{prop2structure1}), (\ref{prop2structureasteriscoprima}) and (\ref{prop2structure2}), we get the stated result for $\underline{d}=(0,0,d)$.

Now, suppose inductively that, for a fixed $0\leq l\leq d-1$, the basis $\beta^{1}_{\underline{d}'''}$ can be chosen as
\begin{equation}\label{prop2structure5}
\beta^{1}_{\underline{d}'''}=\{\varphi_{\underline{d}^{\alpha},\underline{d}'''}(s^{\alpha}_{z})\,/\,\,\underline{d}^{\alpha}\in M^{1}_{\underline{d}'''}\,\,\text{and}\,\,z=1,\ldots,r_{\alpha}\}.
\end{equation}
Since $\underline{d}=(0,d-l,l)$, we have $V_{\underline{d}}^{X_{2},0}\subseteq V_{\underline{d}}^{X_{1},0}$. Then
\begin{equation}\label{prop2structure3}
V_{\underline{d}}^{X_{1},0}+V_{\underline{d}}^{X_{2},0}+V_{\underline{d}}^{X_{3},0}=V_{\underline{d}}^{X_{1},0}+V_{\underline{d}}^{X_{3},0}. 
\end{equation}
Now, as in the proof of Proposition \ref{prop1structure}, we have
\begin{center}
$V_{\underline{d}}^{X_{3},0}=\varphi_{\underline{d}''',\underline{d}}(W^{1}_{\underline{d}'''})\oplus V_{\underline{d}}^{X_{2}^{c},0}$.
\end{center}
Then
\begin{equation}\label{prop2structure4}
V_{\underline{d}}^{X_{1},0}+V_{\underline{d}}^{X_{3},0}=\varphi_{\underline{d}''',\underline{d}}(W^{1}_{\underline{d}'''})\oplus V_{\underline{d}}^{X_{1},0}.
\end{equation}
From (\ref{prop2structure3}) and (\ref{prop2structure4}), we obtain
\begin{center}
$V_{\underline{d}}^{X_{1},0}+V_{\underline{d}}^{X_{2},0}+V_{\underline{d}}^{X_{3},0}=\varphi_{\underline{d}''',\underline{d}}(W^{1}_{\underline{d}'''})\oplus V_{\underline{d}}^{X_{1},0}$.
\end{center}
Then, from the proof of Proposition \ref{like simple basis}, we conclude that $\beta^{1}_{\underline{d}}$ can be chosen as
\begin{equation}\label{prop2structure6}
\beta^{1}_{\underline{d}}=\varphi_{\underline{d}''',\underline{d}}(\beta^{1}_{\underline{d}'''})\cup \{s^{\widetilde{\alpha}}_{1},\ldots,s^{\widetilde{\alpha}}_{r_{\widetilde{\alpha}}}\}\,\,\text{if}\,\,\underline{d}=\underline{d}^{\widetilde{\alpha}}\,\,\text{for some}\,\,\widetilde{\alpha},
\end{equation}
\begin{equation}\label{prop2structure7}
\text{and}\,\,\beta^{1}_{\underline{d}}=\varphi_{\underline{d}''',\underline{d}}(\beta^{1}_{\underline{d}'''})\,\,\text{otherwise}.
\end{equation}
On the other hand, notice that 
\begin{center}
$M^{1}_{\underline{d}}=\{(0,d-\tilde{l},\tilde{l})\in M\,/\,\,l\leq \tilde{l}\}$ and \\
$M^{1}_{\underline{d}'''}=\{(0,d-\tilde{l},\tilde{l})\in M\,/\,\,l+1\leq \tilde{l}\}$,
\end{center}
so
\begin{equation}\label{prop2structure8}
M^{1}_{\underline{d}}=M^{1}_{\underline{d}'''}\cup \{\underline{d}\}.
\end{equation}
From (\ref{prop2structure5}), (\ref{prop2structure6}), (\ref{prop2structure7}) and (\ref{prop2structure8}), we get the stated result for $\underline{d}=(0,d-l,l)$.

Now, suppose inductively that, for a fixed $1\leq i\leq d-1$, the bases $\beta^{1}_{\underline{d}_{0}}$ can be chosen as
\begin{equation}\label{prop2structureasteriscodobleprima}
\beta^{1}_{\underline{d}_{0}}=\{\varphi_{\underline{d}^{\alpha},\underline{d}_{0}}(s^{\alpha}_{z})\,/\,\,\underline{d}^{\alpha}\in M^{1}_{\underline{d}_{0}}\,\,\text{and}\,\,z=1,\ldots,r_{\alpha}\}
\end{equation}
for any $\underline{d}_{0}=(\tilde{i},d-\tilde{i}-\tilde{l},\tilde{l})$ such that $0\leq \tilde{i}\leq i-1$ and $0\leq \tilde{l}\leq d-\tilde{i}$. Then, for the fixed integer $i$, we will prove the stated result inductively as follows:

For $l=d-i$, we have $\underline{d}=(i,0,d-i)$. Then $V_{\underline{d}}^{X_{1},0}\subseteq V_{\underline{d}}^{X_{2},0}$ and $V_{\underline{d}}^{X_{3},0}\subseteq V_{\underline{d}}^{X_{2},0}$. It follows that
\begin{equation}\label{prop2structure9}
V_{\underline{d}}^{X_{1},0}=V_{\underline{d}}^{X_{3}^{c},0}\,\,\text{and}\,\,V_{\underline{d}}^{X_{1},0}+V_{\underline{d}}^{X_{2},0}+V_{\underline{d}}^{X_{3},0}=V_{\underline{d}}^{X_{2},0}.
\end{equation}
Now, as in the proof of Proposition \ref{prop1structure}, we have
\begin{equation}\label{prop2structure10}
V_{\underline{d}}^{X_{2},0}=\varphi_{\underline{d}',\underline{d}}(W^{1}_{\underline{d}'})\oplus V_{\underline{d}}^{X_{3}^{c},0}.
\end{equation}
From (\ref{prop2structure9}) and (\ref{prop2structure10}), we obtain
\begin{center}
$V_{\underline{d}}^{X_{1},0}+V_{\underline{d}}^{X_{2},0}+V_{\underline{d}}^{X_{3},0}=\varphi_{\underline{d}',\underline{d}}(W^{1}_{\underline{d}'})\oplus V_{\underline{d}}^{X_{1},0}$.
\end{center}
Then, from the proof of Proposition \ref{like simple basis}, we conclude that $\beta^{1}_{\underline{d}}$ can be chosen as
\begin{equation}\label{prop2structure11}
\beta^{1}_{\underline{d}}=\varphi_{\underline{d}',\underline{d}}(\beta^{1}_{\underline{d}'})\cup \{s^{\widetilde{\alpha}}_{1},\ldots,s^{\widetilde{\alpha}}_{r_{\widetilde{\alpha}}}\}\,\,\text{if}\,\,\underline{d}=\underline{d}^{\widetilde{\alpha}}\,\,\text{for some}\,\,\widetilde{\alpha},
\end{equation}
\begin{equation}\label{prop2structure12}
\text{and}\,\,\beta^{1}_{\underline{d}}=\varphi_{\underline{d}',\underline{d}}(\beta^{1}_{\underline{d}'})\,\,\text{otherwise}.
\end{equation}
On the other hand, since
\begin{center}
$M^{1}_{\underline{d}}=\{(\tilde{i},d-\tilde{i}-\tilde{l},\tilde{l})\in M\,/\,\,\tilde{i}\leq i\,\,\text{and}\,\,\tilde{i}-i\leq \tilde{l}-(d-i)\}$ and \\
$M^{1}_{\underline{d}'}=\{(\tilde{i},d-\tilde{i}-\tilde{l},\tilde{l})\in M\,/\,\,\tilde{i}\leq i-1\,\,\text{and}\,\,\tilde{i}-(i-1)\leq \tilde{l}-(d-i-1)\}$,
\end{center}
we have
\begin{equation}\label{prop2structure13}
M^{1}_{\underline{d}}=M^{1}_{\underline{d}'}\cup \{\underline{d}\}.
\end{equation}
From the induction hypothesis (\ref{prop2structureasteriscodobleprima}) applied to $\underline{d}_{0}=\underline{d}'$, (\ref{prop2structure11}), (\ref{prop2structure12}) and (\ref{prop2structure13}), we get the stated result for $\underline{d}=(i,0,d-i)$.

Now, suppose inductively that, for a fixed $1\leq l\leq d-i-1$, the basis $\beta^{1}_{\underline{d}'''}$ can be chosen as
\begin{equation}\label{prop2structure14}
\beta^{1}_{\underline{d}'''}=\{\varphi_{\underline{d}^{\alpha},\underline{d}'''}(s^{\alpha}_{z})\,/\,\,\underline{d}^{\alpha}\in M^{1}_{\underline{d}'''}\,\,\text{and}\,\,z=1,\ldots,r_{\alpha}\}.
\end{equation}
As in the proof of Proposition \ref{like simple basis}, we have
\begin{equation}\label{prop2structure15}
V_{\underline{d}}^{X_{2},0}+V_{\underline{d}}^{X_{3},0}=(V_{\underline{d}}^{X_{3}^{c},0}\oplus V_{\underline{d}}^{X_{2}^{c},0})\oplus (\varphi_{\underline{d}',\underline{d}}(W^{1}_{\underline{d}'})+\varphi_{\underline{d}''',\underline{d}}(W^{1}_{\underline{d}'''})).
\end{equation} 
Since the distributivity holds in $V_{\underline{d}}$, we have that
\begin{equation}\label{prop2structure16}
V_{\underline{d}}^{X_{1},0}\cap (V_{\underline{d}}^{X_{2},0}+V_{\underline{d}}^{X_{3},0})=V_{\underline{d}}^{X_{3}^{c},0}\oplus V_{\underline{d}}^{X_{2}^{c},0}.
\end{equation}
Thus, from (\ref{prop2structure15}) and (\ref{prop2structure16}), we obtain
\begin{center}
$V_{\underline{d}}^{X_{1},0}+V_{\underline{d}}^{X_{2},0}+V_{\underline{d}}^{X_{3},0}=(\varphi_{\underline{d}',\underline{d}}(W^{1}_{\underline{d}'})+\varphi_{\underline{d}''',\underline{d}}(W^{1}_{\underline{d}'''}))\oplus V_{\underline{d}}^{X_{1},0}$.
\end{center}
Then, from the proof of Proposition \ref{like simple basis}, we conclude that $\beta^{1}_{\underline{d}}$ can be chosen as
\begin{equation}\label{prop2structure17}
\beta^{1}_{\underline{d}}=(\varphi_{\underline{d}',\underline{d}}(\beta^{1}_{\underline{d}'})\cup \varphi_{\underline{d}''',\underline{d}}(\beta^{1}_{\underline{d}'''}))\cup \{s^{\widetilde{\alpha}}_{1},\ldots,s^{\widetilde{\alpha}}_{r_{\widetilde{\alpha}}}\}\,\,\text{if}\,\,\underline{d}=\underline{d}^{\widetilde{\alpha}}\,\,\text{for some}\,\,\widetilde{\alpha},
\end{equation}
\begin{equation}\label{prop2structure18}
\text{and}\,\,\beta^{1}_{\underline{d}}=\varphi_{\underline{d}',\underline{d}}(\beta^{1}_{\underline{d}'})\cup \varphi_{\underline{d}''',\underline{d}}(\beta^{1}_{\underline{d}'''})\,\,\text{otherwise}.
\end{equation}
On the other hand, it is easily seen that
\begin{equation}\label{prop2structure19}
M^{1}_{\underline{d}}=(M^{1}_{\underline{d}'}\cup M^{1}_{\underline{d}'''})\cup \{\underline{d}\}.
\end{equation}
From the induction hypothesis (\ref{prop2structureasteriscodobleprima}) applied to $\underline{d}_{0}=\underline{d}'$, (\ref{prop2structure14}), (\ref{prop2structure17}), (\ref{prop2structure18}) and (\ref{prop2structure19}), we get the stated result for $\underline{d}=(i,d-i-l,l)$.

As for $l=0$, we have $\underline{d}=(i,d-i,0)$. Then $V_{\underline{d}}^{X_{2},0}\subseteq V_{\underline{d}}^{X_{3},0}$, and hence
\begin{center}
$V_{\underline{d}}^{X_{1},0}+V_{\underline{d}}^{X_{2},0}+V_{\underline{d}}^{X_{3},0}=V_{\underline{d}}^{X_{1},0}+V_{\underline{d}}^{X_{3},0}$. 
\end{center}
From here, the proof follows as in the case of the multidegree $(0,d-l,l)$, for $0\leq l\leq d-1$, as in the present case the result holds for $\underline{d}'''=(i,d-i-1,1)$, and we have $M^{1}_{\underline{d}}=M^{1}_{\underline{d}'''}\cup \{\underline{d}\}$ as well.

Finally, for $i=d$ and $l=0$, by Proposition \ref{prop1structure}, item 6, we have
\begin{center}
$V_{\underline{d}}^{X_{1},0}+V_{\underline{d}}^{X_{2},0}+V_{\underline{d}}^{X_{3},0}=\varphi_{\underline{d}'',\underline{d}}(W^{1}_{\underline{d}''})$.
\end{center}
Then, from the proof of Proposition \ref{like simple basis}, we conclude that $\beta^{1}_{\underline{d}}$ can be chosen as
\begin{equation}\label{prop2structure20}
\beta^{1}_{\underline{d}}=\varphi_{\underline{d}'',\underline{d}}(\beta^{1}_{\underline{d}''})\cup \{s^{\widetilde{\alpha}}_{1},\ldots,s^{\widetilde{\alpha}}_{r_{\widetilde{\alpha}}}\}\,\,\text{if}\,\,\underline{d}=\underline{d}^{\widetilde{\alpha}}\,\,\text{for some}\,\,\widetilde{\alpha},
\end{equation}
\begin{equation}\label{prop2structure21}
\text{and}\,\,\beta^{1}_{\underline{d}}=\varphi_{\underline{d}'',\underline{d}}(\beta^{1}_{\underline{d}''})\,\,\text{otherwise}.
\end{equation}
On the other hand, it is easily seen that
\begin{equation}\label{prop2structure22}
M^{1}_{\underline{d}}=M^{1}_{\underline{d}''}\cup \{\underline{d}\}.
\end{equation}
Since the result holds for $\underline{d}''=(d-1,1,0)$, it follows from (\ref{prop2structure20}), (\ref{prop2structure21}) and (\ref{prop2structure22}) that the result holds for $\underline{d}=(d,0,0)$ as well. This finishes the proof for the case $q=1$.

Now, we will prove the case $q=2$. For $i+l=d$, we have $\underline{d}=(i,0,d-i)$, so $V_{\underline{d}}^{X_{1},0}\subseteq V_{\underline{d}}^{X_{2},0}$ and $V_{\underline{d}}^{X_{3},0}\subseteq V_{\underline{d}}^{X_{2},0}$. It follows that
\begin{center}
$V_{\underline{d}}^{X_{1},0}+V_{\underline{d}}^{X_{2},0}+V_{\underline{d}}^{X_{3},0}=V_{\underline{d}}^{X_{2},0}$.
\end{center}
Then, from the proof of Proposition \ref{like simple basis}, we conclude that $\beta^{2}_{\underline{d}}$ can be chosen as
\begin{equation}\label{prop2structure23}
\beta^{2}_{\underline{d}}=\{s^{\widetilde{\alpha}}_{1},\ldots,s^{\widetilde{\alpha}}_{r_{\widetilde{\alpha}}}\}\,\,\text{if}\,\,\underline{d}=\underline{d}^{\widetilde{\alpha}}\,\,\text{for some}\,\,\widetilde{\alpha},
\end{equation}
\begin{equation}\label{prop2structure24}
\text{and}\,\,\beta^{2}_{\underline{d}}=\emptyset\,\,\text{otherwise}.
\end{equation}
On the other hand, notice that 
\begin{equation}\label{prop2structure25}
M^{2}_{\underline{d}}=\{\underline{d}\}, 
\end{equation}
as $i+l=d$. Thus, if $\underline{d}=\underline{d}^{\widetilde{\alpha}}$ for some $\widetilde{\alpha}$, we have
\begin{equation}\label{prop2structure26repetido}
\{s^{\widetilde{\alpha}}_{1},\ldots,s^{\widetilde{\alpha}}_{r_{\widetilde{\alpha}}}\}=\{\varphi_{\underline{d}^{\alpha},\underline{d}}(s^{\alpha}_{z})\,/\,\,\underline{d}^{\alpha}\in M^{2}_{\underline{d}}\,\,\text{and}\,\,z=1,\ldots,r_{\alpha}\}.
\end{equation}
From (\ref{prop2structure23}), (\ref{prop2structure24}), (\ref{prop2structure25}) and (\ref{prop2structure26repetido}), we get the stated result for $\underline{d}=(i,0,d-i)$.

Now, suppose inductively that, for a fixed $1\leq m\leq d$, the result holds for any multidegree $(\tilde{i},d-\tilde{i}-\tilde{l},\tilde{l})$ such that $\tilde{i}+\tilde{l}\geq d-(m-1)$. Then, we will prove the stated result for $\underline{d}=(i,d-i-l,l)$, with $i+l=d-m$.

Set $\widetilde{\underline{d}}:=(i+1,d-i-l-1,l)$. It follows from the proof of Proposition \ref{like simple basis} that
\begin{equation}\label{prop2structure26}
V_{\underline{d}}^{X_{1},0}+V_{\underline{d}}^{X_{3},0}=(V_{\underline{d}}^{X_{3}^{c},0}\oplus V_{\underline{d}}^{X_{1}^{c},0})\oplus (\varphi_{\widetilde{\underline{d}},\underline{d}}(W^{2}_{\widetilde{\underline{d}}})+\varphi_{\underline{d}''',\underline{d}}(W^{2}_{\underline{d}'''})).
\end{equation} 
Since the distributivity holds in $V_{\underline{d}}$, we have that
\begin{equation}\label{prop2structure27}
V_{\underline{d}}^{X_{2},0}\cap (V_{\underline{d}}^{X_{1},0}+V_{\underline{d}}^{X_{3},0})=V_{\underline{d}}^{X_{3}^{c},0}\oplus V_{\underline{d}}^{X_{1}^{c},0}.
\end{equation}
Thus, from (\ref{prop2structure26}) and (\ref{prop2structure27}), we obtain
\begin{center}
$V_{\underline{d}}^{X_{1},0}+V_{\underline{d}}^{X_{2},0}+V_{\underline{d}}^{X_{3},0}=(\varphi_{\widetilde{\underline{d}},\underline{d}}(W^{2}_{\widetilde{\underline{d}}})+\varphi_{\underline{d}''',\underline{d}}(W^{2}_{\underline{d}'''}))\oplus V_{\underline{d}}^{X_{2},0}$.
\end{center}
Then, from the proof of Proposition \ref{like simple basis}, we conclude that $\beta^{2}_{\underline{d}}$ can be chosen as
\begin{equation}\label{prop2structure28}
\beta^{2}_{\underline{d}}=(\varphi_{\widetilde{\underline{d}},\underline{d}}(\beta^{2}_{\widetilde{\underline{d}}})\cup \varphi_{\underline{d}''',\underline{d}}(\beta^{2}_{\underline{d}'''}))\cup \{s^{\widetilde{\alpha}}_{1},\ldots,s^{\widetilde{\alpha}}_{r_{\widetilde{\alpha}}}\}\,\,\text{if}\,\,\underline{d}=\underline{d}^{\widetilde{\alpha}}\,\,\text{for some}\,\,\widetilde{\alpha},
\end{equation}
\begin{equation}\label{prop2structure29}
\text{and}\,\,\beta^{2}_{\underline{d}}=\varphi_{\widetilde{\underline{d}},\underline{d}}(\beta^{2}_{\widetilde{\underline{d}}})\cup \varphi_{\underline{d}''',\underline{d}}(\beta^{2}_{\underline{d}'''})\,\,\text{otherwise}.
\end{equation}
On the other hand, it is easily seen that
\begin{equation}\label{prop2structure30}
M^{2}_{\underline{d}}=(M^{2}_{\widetilde{\underline{d}}}\cup M^{2}_{\underline{d}'''})\cup \{\underline{d}\}.
\end{equation}
From the induction hypothesis applied to $\widetilde{\underline{d}}$ and $\underline{d}'''$, (\ref{prop2structure28}), (\ref{prop2structure29}) and (\ref{prop2structure30}), we get the stated result for $\underline{d}=(i,d-i-l,l)$. This finishes the proof for the case $q=2$.
\hfill $\Box$

Using the three previous results, we obtain the following proposition.

\begin{prop}\label{simple basis existence}
Let $\{(\mathcal{L}_{\underline{d}},V_{\underline{d}})\}_{\underline{d}}$ be an exact limit linear series, and suppose that the distributivity holds in $V_{\underline{d}}$ for any $\underline{d}$. Then $\{(\mathcal{L}_{\underline{d}},V_{\underline{d}})\}_{\underline{d}}$ is simple.
\end{prop}
{\em Proof.} Write
\begin{center}
$\{\underline{d}\,/\,\,V_{\underline{d}}^{X_{1},0}+V_{\underline{d}}^{X_{2},0}+V_{\underline{d}}^{X_{3},0}\neq V_{\underline{d}}\}=\{\underline{d}^{1},\ldots,\underline{d}^{m}\}$,
\end{center}
and for each $\alpha=1,\ldots,m$, let $\{s^{\alpha}_{1},\ldots,s^{\alpha}_{r_{\alpha}}\}\subseteq V_{\underline{d}^{\alpha}}$ be any linearly independent subset such that
\begin{center}
$V_{\underline{d}^{\alpha}}=(V_{\underline{d}}^{X_{1},0}+V_{\underline{d}}^{X_{2},0}+V_{\underline{d}}^{X_{3},0})\oplus \left\langle s^{\alpha}_{1},\ldots,s^{\alpha}_{r_{\alpha}}\right\rangle$.
\end{center}
It follows from Corollary \ref{sum of codimensions} that
\begin{equation}\label{simple basis existence1}
r_{1}+\ldots+r_{m}=r+1.
\end{equation}
Now, we claim that
\begin{equation}\label{simple basis existence2}
V_{\underline{d}}=\left\langle \varphi_{\underline{d}^{1},\underline{d}}(s^{1}_{1}),\ldots,\varphi_{\underline{d}^{1},\underline{d}}(s^{1}_{r_{1}}),\ldots,\varphi_{\underline{d}^{m},\underline{d}}(s^{m}_{1}),\ldots,\varphi_{\underline{d}^{m},\underline{d}}(s^{m}_{r_{m}})\right\rangle,
\end{equation}
for any $\underline{d}=(i,d-i-l,l)$. In fact, as in the statement of Proposition \ref{prop1structure}, there are six cases to consider. We will only prove the claim for $i\geq 1$, $l\geq 1$ and $i+l\leq d-1$, as the remaining cases are analogous. In this case, by Proposition \ref{prop1structure}, item 4, we get
\begin{align}\label{simple basis existence asterisco}
V_{\underline{d}}^{X_{1},0}+V_{\underline{d}}^{X_{2},0}+V_{\underline{d}}^{X_{3},0}= & \,\varphi_{\underline{d}',\underline{d}}(W^{1}_{\underline{d}'})+\varphi_{\underline{d}',\underline{d}}(W^{3}_{\underline{d}'}) \nonumber \\
& +\varphi_{\underline{d}''',\underline{d}}(W^{1}_{\underline{d}'''})+\varphi_{\underline{d}''',\underline{d}}(W^{2}_{\underline{d}'''})+\varphi_{\widetilde{\underline{d}},\underline{d}}(W^{2}_{\widetilde{\underline{d}}})+ \varphi_{\widetilde{\underline{d}},\underline{d}}(W^{3}_{\widetilde{\underline{d}}}),
\end{align} 
where $\widetilde{\underline{d}}:=(i+1,d-i-l-1,l)$. Also, it is easily seen that
\begin{center}
$M^{1}_{\underline{d}}=(M^{1}_{\underline{d}'}\cup M^{1}_{\underline{d}'''})\cup \{\underline{d}\}$, \\
$M^{2}_{\underline{d}}=(M^{2}_{\widetilde{\underline{d}}}\cup M^{2}_{\underline{d}'''})\cup \{\underline{d}\}$ \\
and $M^{3}_{\underline{d}}=(M^{3}_{\underline{d}'}\cup M^{3}_{\widetilde{\underline{d}}})\cup \{\underline{d}\}$,
\end{center}
that is,
\begin{equation}\label{simple basis existence3}
M^{1}_{\underline{d}}\setminus\{\underline{d}\}=M^{1}_{\underline{d}'}\cup M^{1}_{\underline{d}'''},
\end{equation}
\begin{equation}\label{simple basis existence4}
M^{2}_{\underline{d}}\setminus\{\underline{d}\}=M^{2}_{\widetilde{\underline{d}}}\cup M^{2}_{\underline{d}'''}
\end{equation}
\begin{equation}\label{simple basis existence5}
\text{and}\,\,M^{3}_{\underline{d}}\setminus\{\underline{d}\}=M^{3}_{\underline{d}'}\cup M^{3}_{\widetilde{\underline{d}}}.
\end{equation} 
By Proposition \ref{prop2structure}, the bases $\beta^{1}_{\underline{d}'}$ and $\beta^{1}_{\underline{d}'''}$ can be chosen as
\begin{equation}\label{simple basis existence6}
\beta^{1}_{\underline{d}'}=\{\varphi_{\underline{d}^{\alpha},\underline{d}'}(s^{\alpha}_{z})\,/\,\,\underline{d}^{\alpha}\in M^{1}_{\underline{d}'}\,\,\text{and}\,\,z=1,\ldots,r_{\alpha}\}
\end{equation}
\begin{equation}\label{simple basis existence7}
\text{and}\,\,\beta^{1}_{\underline{d}'''}=\{\varphi_{\underline{d}^{\alpha},\underline{d}'''}(s^{\alpha}_{z})\,/\,\,\underline{d}^{\alpha}\in M^{1}_{\underline{d}'''}\,\,\text{and}\,\,z=1,\ldots,r_{\alpha}\}.
\end{equation}
From (\ref{simple basis existence3}), (\ref{simple basis existence6}) and (\ref{simple basis existence7}), we conclude that
\begin{equation}\label{simple basis existence8}
\varphi_{\underline{d}',\underline{d}}(W^{1}_{\underline{d}'})+\varphi_{\underline{d}''',\underline{d}}(W^{1}_{\underline{d}'''})=\left\langle\{\varphi_{\underline{d}^{\alpha},\underline{d}}(s^{\alpha}_{z})\,/\,\,\underline{d}^{\alpha}\in M^{1}_{\underline{d}}\setminus\{\underline{d}\}\,\,\text{and}\,\,z=1,\ldots,r_{\alpha}\}\right\rangle.
\end{equation}
Analogously, using Proposition \ref{prop2structure}, (\ref{simple basis existence4}) and (\ref{simple basis existence5}), we get
\begin{equation}\label{simple basis existence9}
\varphi_{\underline{d}''',\underline{d}}(W^{2}_{\underline{d}'''})+\varphi_{\widetilde{\underline{d}},\underline{d}}(W^{2}_{\widetilde{\underline{d}}})=\left\langle\{\varphi_{\underline{d}^{\alpha},\underline{d}}(s^{\alpha}_{z})\,/\,\,\underline{d}^{\alpha}\in M^{2}_{\underline{d}}\setminus\{\underline{d}\}\,\,\text{and}\,\,z=1,\ldots,r_{\alpha}\}\right\rangle
\end{equation}
and
\begin{equation}\label{simple basis existence10}
\varphi_{\underline{d}',\underline{d}}(W^{3}_{\underline{d}'})+ \varphi_{\widetilde{\underline{d}},\underline{d}}(W^{3}_{\widetilde{\underline{d}}})=\left\langle\{\varphi_{\underline{d}^{\alpha},\underline{d}}(s^{\alpha}_{z})\,/\,\,\underline{d}^{\alpha}\in M^{3}_{\underline{d}}\setminus\{\underline{d}\}\,\,\text{and}\,\,z=1,\ldots,r_{\alpha}\}\right\rangle.
\end{equation}
Since $M^{1}_{\underline{d}}\cup M^{2}_{\underline{d}}\cup M^{3}_{\underline{d}}=M$, it follows from (\ref{simple basis existence asterisco}), (\ref{simple basis existence8}), (\ref{simple basis existence9}) and (\ref{simple basis existence10}) that
\begin{equation}\label{simple basis existence11}
V_{\underline{d}}^{X_{1},0}+V_{\underline{d}}^{X_{2},0}+V_{\underline{d}}^{X_{3},0}=\left\langle\{\varphi_{\underline{d}^{\alpha},\underline{d}}(s^{\alpha}_{z})\,/\,\,\underline{d}^{\alpha}\in M\setminus\{\underline{d}\}\,\,\text{and}\,\,z=1,\ldots,r_{\alpha}\}\right\rangle.
\end{equation} 
If $\underline{d}=\underline{d}^{\widetilde{\alpha}}$ for some $\widetilde{\alpha}$, we have
\begin{equation}\label{simple basis existence12}
V_{\underline{d}}=(V_{\underline{d}}^{X_{1},0}+V_{\underline{d}}^{X_{2},0}+V_{\underline{d}}^{X_{3},0})\oplus \left\langle s^{\widetilde{\alpha}}_{1},\ldots,s^{\widetilde{\alpha}}_{r_{\widetilde{\alpha}}}\right\rangle.
\end{equation}
Then, from (\ref{simple basis existence11}) and (\ref{simple basis existence12}), we get
\begin{center}
$V_{\underline{d}}=\left\langle\{\varphi_{\underline{d}^{\alpha},\underline{d}}(s^{\alpha}_{z})\,/\,\,\underline{d}^{\alpha}\in M\,\,\text{and}\,\,z=1,\ldots,r_{\alpha}\}\right\rangle$.
\end{center}
Now, if $\underline{d}\neq \underline{d}^{\alpha}$ for any $\alpha$, it follows from (\ref{simple basis existence11}) that
\begin{center}
$V_{\underline{d}}^{X_{1},0}+V_{\underline{d}}^{X_{2},0}+V_{\underline{d}}^{X_{3},0}=\left\langle\{\varphi_{\underline{d}^{\alpha},\underline{d}}(s^{\alpha}_{z})\,/\,\,\underline{d}^{\alpha}\in M\,\,\text{and}\,\,z=1,\ldots,r_{\alpha}\}\right\rangle$.
\end{center}
Since $V_{\underline{d}}=V_{\underline{d}}^{X_{1},0}+V_{\underline{d}}^{X_{2},0}+V_{\underline{d}}^{X_{3},0}$, as $\underline{d}\neq \underline{d}^{\alpha}$ for any $\alpha$, we conclude that
\begin{center}
$V_{\underline{d}}=\left\langle\{\varphi_{\underline{d}^{\alpha},\underline{d}}(s^{\alpha}_{z})\,/\,\,\underline{d}^{\alpha}\in M\,\,\text{and}\,\,z=1,\ldots,r_{\alpha}\}\right\rangle$.
\end{center}
So our claim is established.

Since $\text{dim}\,V_{\underline{d}}=r+1$, it follows from (\ref{simple basis existence1}) and (\ref{simple basis existence2}) that 
\begin{center}
$\varphi_{\underline{d}^{1},\underline{d}}(s^{1}_{1}),\ldots,\varphi_{\underline{d}^{1},\underline{d}}(s^{1}_{r_{1}}),\ldots,\varphi_{\underline{d}^{m},\underline{d}}(s^{m}_{1}),\ldots,\varphi_{\underline{d}^{m},\underline{d}}(s^{m}_{r_{m}})$ 
\end{center}
form a basis for $V_{\underline{d}}$, for each $\underline{d}$. Then $\{(\mathcal{L}_{\underline{d}},V_{\underline{d}})\}_{\underline{d}}$ is simple.
\hfill $\Box$

Finally, we obtain the following characterization of simple limit linear series among the exact limit linear series.

\begin{cor}\label{simple characterization}
Let $\{(\mathcal{L}_{\underline{d}},V_{\underline{d}})\}_{\underline{d}}$ be an exact limit linear series of degree $d$ and dimension $r$. Then the following statements are equivalent:\\
1. $\{(\mathcal{L}_{\underline{d}},V_{\underline{d}})\}_{\underline{d}}$ is simple.\\
2. {\em $\displaystyle\sum_{\underline{d}}\text{dim}\,\left(\dfrac{V_{\underline{d}}}{V_{\underline{d}}^{X_{1},0}+V_{\underline{d}}^{X_{2},0}+V_{\underline{d}}^{X_{3},0}}\right)=r+1$}.\\
3. The distributivity holds in $V_{\underline{d}}$, for any $\underline{d}$.
\end{cor}
{\em Proof.} This follows immediately from Proposition \ref{multidegrees characterization}, Corollary \ref{sum of codimensions} and Proposition \ref{simple basis existence}.
\hfill $\Box$

\vspace{0.3cm}

\textbf{Disclosure statement.} The author reports there are no competing interests to declare.

\bibliographystyle{alpha}

\end{document}